\documentclass[11pt]{amsart}
\usepackage{amsxtra}
\usepackage{amssymb}
\usepackage{mathtools}
\theoremstyle{plain}

\newcommand{\cleqn}{\setcounter{equation}{0}}
\newcommand{\clth}{\setcounter{theorem}{0}}
\newcommand {\sectionnew}[1]{\section{#1}\cleqn\clth}

\newtheorem{theorem}{Theorem}[section]
\newtheorem{lemma}[theorem]{Lemma}
\newtheorem{definition-theorem}[theorem]{Definition-Theorem}
\newtheorem{proposition}[theorem]{Proposition}
\newtheorem{problem}[theorem]{Problem}
\newtheorem{corollary}[theorem]{Corollary}
\newtheorem{definition}[theorem]{Definition}
\newtheorem{example}[theorem]{Example}
\newtheorem{remark}[theorem]{Remark}
\newtheorem{conjecture}[theorem]{Conjecture}

\newcommand \bth[1] { \begin{theorem}\label{t#1} }
\newcommand \ble[1] { \begin{lemma}\label{l#1} }

\newcommand \bpr[1] { \begin{proposition}\label{p#1} }
\newcommand \bprob[1] { \begin{problem}\label{q#1} }
\newcommand \bco[1] { \begin{corollary}\label{c#1} }
\newcommand \bde[1] { \begin{definition}\label{d#1}\rm }
\newcommand \bex[1] { \begin{example}\label{e#1}\rm }
\newcommand \bre[1] { \begin{remark}\label{r#1}\rm }
\newcommand \bcj[1] { \begin{conjecture}\label{j#1}\rm }
\renewcommand {\eth} { \end{theorem} }
\newcommand {\ele} { \end{lemma} }

\newcommand {\epr} { \end{proposition} }
\newcommand {\eprob} {\end{problem} }
\newcommand {\eco} { \end{corollary} }
\newcommand {\ede} { \end{definition} }
\newcommand {\eex} { \end{example} }
\newcommand {\ere} { \end{remark} }
\newcommand {\ecj} { \end{conjecture} }

\newcommand {\enota} { \end{notation} }
\newcommand \thref[1]{Theorem \ref{t#1}}
\newcommand \leref[1]{Lemma \ref{l#1}}
\newcommand \prref[1]{Proposition \ref{p#1}}
\newcommand \probref[1]{Problem \ref{q#1}}
\newcommand \coref[1]{Corollary \ref{c#1}}

\newcommand \reref[1]{Remark \ref{r#1}}
\newcommand \lb[1]{\label{#1}}

\def \Rset {{\mathbb R}}         
\def \Cset {{\mathbb C}}
\def \Aset {{\mathbb A}}
\def \KK {{\mathbb K}}
\def \Zset {{\mathbb Z}}
\def \Nset {{\mathbb N}}

\def \AA  {{\mathcal{A}}}           

\def \CC {\mathcal{C}}
\def \FF {{\mathcal{F}}}
\def \II {{\mathcal{I}}}

\def \PP {{\mathcal{P}}}

\def \UU {{\mathcal{U}}}

\def \ZZ {{\mathcal{Z}}}

\def \ZZ {{\mathcal{Z}}}

\def \De {\Delta}   
\def \de {\delta}
\def \al {\alpha}
\def \be {\beta}
\def \vpi {\varpi}
\def \la {\lambda}
\def \La {\Lambda}

\def \Om {\Omega}
\def \ga {\gamma}
\def \de {\delta}
\def \Ga {\Gamma}
\def \Sig {\Sigma}

\def \mt  {\mapsto}
\def \lra {\longrightarrow}

\def \hra {\hookrightarrow}

\def \ci  {\circ}

\def \rcor {\rangle}
\def \lcor {\langle}

\def \ol {\overline}
\def \wt {\widetilde}
\def \wh {\widehat}

\def \id { {\mathrm{id}} }

\def \Ad { {\mathrm{Ad}} }

\def \rank { {\mathrm{rank}} }
\def \Lie { {\mathrm{Lie \,}} }

\def \g  {\mathfrak{g}}   

\def \sl {\mathfrak{sl}} 
\def \h  {\mathfrak{h}}

\def \n  {\mathfrak{n}}

\def \b  {\mathfrak{b}}

\def \sl {\mathfrak{sl}}

\DeclareMathOperator \Span { {\mathrm{Span}} }
\DeclareMathOperator \Aut { {\mathrm{Aut}} }
\DeclareMathOperator \UAut { {\mathrm{UAut}}}
\DeclareMathOperator  \MAut { {\mathrm{Mod-Aut}}}

\DeclareMathOperator \rad { {\mathrm{rad}}} 
\DeclareMathOperator \htt { {\mathrm{ht}}}

\DeclareMathOperator \opp { {\mathrm{op}} }

\DeclareMathOperator \diag { {\mathrm{diag}} }
\DeclareMathOperator \Ker { {\mathrm{Ker}} }

\DeclareMathOperator \Con { {\mathrm{Cone}}}

\DeclareMathOperator \st   { {\mathrm{st}} }

\DeclareMathOperator \Supp { {\mathrm{Supp}} }

\renewcommand \max { {\mathrm{max}} }
\newcommand \Spec { {\mathrm{Spec}} }
\begin{document}
\title[Rigidity of Quadratic Poisson tori]
{Rigidity of Quadratic Poisson tori}
\dedicatory{To George Lusztig on his 70th birthday, with admiration}
\author[Jesse Levitt]{Jesse Levitt}
\address{
Department of Mathematics \\
University of Southern California, 3620 S. Vermont Ave \\
KAP 108, Los Angeles, CA 90089-2532, USA
}
\email{jslevitt@usc.edu}
\author[Milen Yakimov]{Milen Yakimov}
\address{
Department of Mathematics \\
Louisiana State University \\
Baton Rouge, LA 70803,
USA
}
\email{yakimov@math.lsu.edu}
\thanks{The research of J.L. has been supported by a GAANN fellowship and a VIGRE fellowship through the NSF grant DMS-0739382 
and that of M.Y. by the NSF grant DMS-1601862.}
\keywords{Automorphisms of Poisson algebras, cluster algebras, quadratic Poisson tori, Schubert cells}
\subjclass[2000]{Primary: 53D17; Secondary: 16W20, 13F60, 17B37}
\begin{abstract} We prove a rigidity theorem for the Poisson automorphisms of the function fields of tori with quadratic 
Poisson structures over fields of characteristic 0. 
It gives an effective method for classifying the full Poisson automorphism groups of $\Nset$-graded connected cluster algebras
equipped with Gekhtman--Shapiro--Vainshtein Poisson structures. Based on this, we classify the groups of algebraic Poisson automorphisms 
of the open Schubert cells of the full flag varieties of semisimple algebraic groups over fields of characteristic 0, equipped with the standard Poisson structures. 
Their coordinate rings can be identified with the semiclassical limits of the positive parts $U_q(\n_+)$ of the quantized universal enveloping algebras 
of semisimple Lie algebras, 
and the last result establishes a Poisson version of the Andruskiewitsch--Dumas conjecture on $\Aut U_q(\n_+)$.
\end{abstract}
\maketitle
\sectionnew{Introduction}
\lb{Intro}
\noindent
\subsection{Automorphism groups of algebras}
\label{1.1}
The automorphism groups of infinite dimensional (commutative and noncommutative) algebras are difficult to describe and are known 
in very few situations. Before the 70's, the automorphisms of the polynomial and free algebras in two generators, as well as the first Weyl algebra, 
were  determined. However, the majority of the subsequent results were on the existence of wild automorphisms: Joseph \cite{J0} for the universal 
enveloping algebra of $\sl_2$, Umirbaev--Shestakov \cite{SU} for the polynomial algebra in three variables, and others. The automorphism groups 
of these algebras and their natural generalizations remain unknown. 

In the early 2000's Andruskiewitsch and Dumas \cite{AD} posed a conjecture about a concrete description of the automorphism groups of the positive 
parts $U_q(\n_+)$ of the quantized universal enveloping algebras $U_q(\g)$ of split simple Lie algebras $\g$ when $q$ is a non-root of unity. 
This is a drastically different behavior
from the case of universal enveloping algebras and polynomial algebras (in dimension $\geq 3$) 
whose automorphism groups are very large and not classifiable with the current 
methods. This and related problems were settled in \cite{Y-LL,Y-AD}, where a general classification method was given for the automorphisms 
of $\Nset$-graded connected quantum cluster algebras. In the root of unity case (i.e., PI algebras), Ceken, Palmieri, Wang, and Zhang \cite{CPWZ1} 
classified the automorphisms of skew-polynomial algebras using discriminants. However, the Lie theoretic side, e.g. $\Aut U_q(\n_+)$ for $\g$ 
of rank at least 2 and $q$ a root of unity remains unknown.
\subsection{A rigidity theorem and automorphisms of Poisson cluster algebras}
\label{1.2}
In this paper we prove a general rigidity theorem for quadratic Poisson structures on tori over fields of characteristic 0. 
We design a classification scheme, based on it, for the Poisson automorphism groups of $\Nset$-graded connected cluster algebras, 
equipped with Gekhtman--Shapiro--Vainshtein Poisson structures. Many important cluster algebras coming from Lie theory fall in this class, 
see the next subsection for details.

Cluster algebras form an axiomatic class of algebras, introduced by Fomin and Zelevinsky \cite{FZ1}.
They are described as the algebras generated by families of polynomial subalgebras obtained from each other by an infinite process 
of mutation, see the expositions \cite{GSVb,W}.
Each cluster $(y_1, \ldots, y_m)$ of such an algebra $\AA$ over a field $\KK$ gives rise to embeddings 
\begin{equation}
\label{emb0}
\KK[y_1, \ldots, y_m] \subset \AA \subset \KK[y_1^{\pm1}, \ldots, y_m^{\pm 1}]
\end{equation}
where the second embedding follows from the Laurent phenomenon \cite{FZ-L} of Fomin and Zelevinsky. A
Gekhtman--Shapiro--Vainshtein (GSV) Poisson structure \cite{GSV0} on $\AA$ is a Poisson structure $\{.,.\}$ which is quadratic in each cluster; 
that is  
\begin{equation}
\label{quad-br}
\{ y_j, y_k \} = \la_{jk} y_j y_k \quad \forall \, 1 \leq k, j \leq m
\end{equation}
for some integral skewsymmetric matrix $(\la_{jk})$. The first and third algebras in \eqref{emb0}
become Poisson algebras under the natural restriction and extension of $\{.,.\}$. 
Denote the $\KK$-torus $T= \Spec \KK[y_1^{\pm 1}, \ldots, y_m^{\pm 1}]$. 
Its rational function field $\KK(T)$ is a Poisson field under the canonical
extension of the Poisson bracket $\{.,.\}$ which will be denoted by the same symbol.

Now assume that a cluster algebra $\AA$ is $\Nset$-graded, $\AA = \oplus_{n \in \Nset} \AA_n$,
and the elements $y_1, \ldots, y_m$ of one of its clusters are homogeneous of positive degree.
Here and below $\Nset := \{0, 1, \ldots\}$. When $\dim \AA_0$ is small (in particular when $\AA$ 
is connected, i.e., $\dim \AA_0 =1$), the description of the automorphism group $\Aut (\AA, \{.,.\})$
of the Poisson algebra $(\AA, \{.,.\})$ boils down to the 
description of the subgroup of unipotent automorphisms, defined by
\[
\UAut( \AA, \{.,.\}) = 
\{ \psi \in \Aut (\AA, \{.,.\}) \mid \psi (f) - f 
\in  \oplus_{n >k} \AA_n, \; \; \forall k \in \Nset, f \in \AA_k \}.
\]
First we prove that every unipotent automorphism $\psi \in \UAut(\AA, \{.,.\})$ has a canonical extension 
$\phi$ to an automorphism of the Poisson function field $(\KK(T), \{.,.\})$ which is unipotent and bi-integral in the 
sense that it satisfies
\begin{itemize}
\item $\psi(y_k), \psi^{-1}(y_k) \in \KK[y_1, \ldots, y_m]$, $1 \leq k \leq m$ (bi-integrality condition), and
\item $\psi(y_k) - y_k \in \KK[T]_{> \deg y_k}$, $1 \leq k \leq m$ (unipotent condition).
\end{itemize}
This embedding result is proved in \prref{aut-emb}. The main result of the paper is the following rigidity of 
unipotent, bi-integral automorphisms of quadratic Poisson tori.
\smallskip
\\
\noindent
{\bf{Theorem A.}} {\em{Let $\KK$ be a field of characteristic 0 and $T$ be a $\KK$-torus, whose coordinate ring 
$\KK[T] = \KK[y_1^{\pm 1}, \ldots, y_m^{\pm 1}]$ is equipped with a quadratic Poisson structure  $\{.,.\}$ 
as in \eqref{quad-br} and a $\Zset$-grading such that $\deg y_k >0$, for all $k$.

Every unipotent, bi-integral automorphism $\phi$ of the Poisson function field $(\KK(T), \{.,.\})$ satisfies
\[
\phi(y_k) y_k^{-1} \in \ZZ(\KK[T]), \quad \forall 1 \leq k \leq m,
\]
where $\ZZ(\KK(T))$ denotes the Poisson center of $(\KK(T), \{.,.\})$.  
}}
\\

This is a rigidity result in the following sense. A quadratic Poisson structure on a torus $T$ 
arises as the semiclassical structure of the deformation of the coordinate ring
$\KK[T]$ to a quantum torus. The theorem proves that the unipotent, bi-integral automorphisms
of the Poisson function field of such a torus only come from its Poisson center which is 
precisely the subalgebra of $\KK[T]$ that is deformed to the center of the corresponding quantum torus.
The theorem is proved in Section \ref{Rig}.

The Poisson automorphisms of $\Nset$-graded cluster algebras $\AA$ with Gekhtman--Shapiro--Vainshtein Poisson structures
$\{.,.\}$ can be classified as follows. We embed $\UAut(\AA, \{.,.\})$ in the set of  
unipotent, bi-integral automorphisms of the Poisson function field of a torus $T$ 
associated to a cluster of $\AA$, and then use the rigidity theorem. In concrete situations $\dim \ZZ(\KK[T])$ is 
much smaller than that of $T$, and one can fully determine $\UAut(\AA, \{.,.\})$ from the action of the 
automorphisms on the height one prime ideals of $\AA$ that are Poisson ideals.
Finally, $\Aut(\AA, \{.,.\})$ is reconstructed from $\UAut(\AA, \{.,.\})$ from the action of 
$\phi \in \Aut(\AA, \{.,.\})$ on $\AA_0$. Section \ref{App-cl} contains 
full details on the method and Section \ref{App-Sc} contains a concrete realization. 
Section \ref{App-cl} also contains further relations to the modular automorphism group 
of the cluster algebra $\AA$ and the group of toric transformations of $\AA$. 
\subsection{Automorphisms of Poisson cluster algebras on unipotent groups}
\label{1.3}
We apply Theorem A and the method of the previous subsection to classify the Poisson automorphisms 
of an important family of cluster algebras on unipotent groups related to dual canonical bases 
and total positivity on flag varieties.

In \cite{L-pos1} Lusztig defined a general notion of total positivity in complex semisimple algebraic groups $G$. 
Fomin and Zelevinsky \cite{FZ0} proved that the nonnegative subset $G_{\geq 0}$  has a canonical decomposition obtained by intersecting $G_{\geq 0}$
with the double Bruhat cells
\[
G^{w,u} := B_+ w B_+ \cap B_- u B_-
\]
where $B_\pm$ is a pair of opposite Borel subgroups and $(w, u)$ is a pair of Weyl group elements.
Berenstein, Fomin and Zelevinsky \cite{BFZ} constructed upper quantum cluster algebra structures $\UU(G^{w,u})$
on the coordinate rings of $G^{w,u}$ that provide a bridge between Poisson structures and total positivity in the following sense:
\begin{enumerate}
\item The set  $G_{\geq 0} \cap G^{w,u}$ is precisely the totally positive part of $G^{w,u}$ in the cluster theoretic sense, i.e., it 
consists of the points on which all cluster variables (equivalently, those in one cluster) take positive values.
\item The double Bruhat cells $G^{w,u}$ are Poisson submanifolds of $G$, equipped with the 
standard (Sklyanin) Poisson structure $\Pi_{\st}$, and the induced Poisson structures on $\Cset[G^{w,u}]$ 
are GSV Poisson structures for $\UU(G^{w,u})$. 
\end{enumerate}

The push forward of $\Pi_{\st}$ under the projection $G \to G/B_+$ is the standard Poisson structure $\pi_{\st}$ on the 
full flag variety $G/B_+$. We restrict $\pi_{\st}$ to the open Schubert cell $B_+ w_\circ B_+/B_+$, identified with the 
unipotent radical $U_+$ of $B_+$ in the usual way, where $w_\circ$ is the longest element of the Weyl group of $G$. Denote the 
corresponding Poisson structures by $\pi_{\st}$:
\[
(U_+, \pi_{\st}) \cong (B_+ w_\circ B_+/B_+, \pi_{\st}) \hra (G/B_+, \pi_{\st}) \twoheadleftarrow (G, \Pi_{\st}).
\]

For a simply laced group $G$, Geiss--Leclerc--Schr\"oer \cite{GLS0} proved that the upper cluster algebra $\UU(G^{w_\circ, 1} )$ 
coincides with the corresponding cluster algebra, and when its frozen variables are not inverted, it gives rise to a cluster algebra
structure $\AA(U_+)$ on $\Cset[U_+]$. (There is a minor detail that for the last relation 
one passes to the reduced double Bruhat cell $G^{w_\circ, 1}/H$ where $H:= B_+ \cap B_-$.) 
For a general semisimple group $G$, these facts were proved in \cite{GY-P}. The cluster algebras 
$\AA(U_+)$ have the following properties:
\begin{enumerate}
\item Geiss, Leclrec and Schr\"oer \cite{GLS1} proved that the cluster monomials of $\AA(U_+)$ belong to Lusztig's dual canonical basis 
\cite{L-semican0,L-semican} of $\Cset[U_+]$ when $G$ is simply laced.
\item The part of the totally nonnegative subset $(G/B_+)_{\geq 0}$ defined and studied by Lusztig \cite{L-pos2}, that is
inside the Schubert cell $B_+ w_\circ B_+/B_+$, is precisely the nonnegative subset in the cluster theoretic 
sense for the cluster algebra $\AA(U_+)$. 
\item The standard Poisson structure $\pi_{\st}$ on $U_+$ is a GSV Poisson structure for $\AA(U_+)$. 
\end{enumerate}
Denote by $\{.,.\}_{\st}$ the Poisson bracket on $\Cset[U_+]$ associated to $\pi_{\st}$. The last result of the paper classifies the Poisson automorphism groups
of the cluster algebras $(\Cset[U_+], \{.,.\}_{\st})$.  
\smallskip
\\
\noindent
{\bf{Theorem B.}}
{\em{Let $G$ be a complex, connected, simply connected, semisimple algebraic group 
which does not have $SL_2$ factors. The automorphism group of the Poisson algebra 
$(\Cset[U_+], \{.,.\}_{\st})$ is isomorphic to 
\[
(H/Z_G) \ltimes \Aut(\Ga),
\]
where $H = B_+ \cap B_-$, $Z_G$ is the center of $G$, $\Ga$ is the Dynkin graph of $G$, 
and $\Aut(\Ga)$ is its automorphism group.
}}
\smallskip
\\

Here $H$ acts on $U_+$ by conjugation and $\Aut(\Ga)$ acts by permuting a fixed set of Chevalley generators. The case when 
$G$ has $SL_2$ factors is excluded because in that case there are pathological problems coming from the fact that 
$(\Cset[U_+], \{.,.\}_{\st})$ is a tensor product of two Poisson algebras, one of which is a polynomial 
algebra with a trivial Poisson bracket. The theorem is proved in Section \ref{App-Sc}. We prove a more general form of the 
theorem for split, connected, simply connected, semisimple algebraic groups $G$ without $SL_2$ factors 
over arbitrary fields $\KK$ of characteristic 0.
It can be viewed as a Poisson analog of the Andruskiewitsch--Dumas conjecture \cite{AD} on $\Aut U_q(\n_+)$. 

Finally, in \thref{5c} we also solve the isomorphism problem for the family of Poisson algebras of the form 
$(\KK[U_+], \{.,.\}_{\st})$ for split, connected, simply connected, semisimple algebraic groups $G$ over 
fields $\KK$ of characteristic 0 (allowing $SL_2$ factors).

We will use the following notation. For $j \leq k \in \Zset$, denote $[j,k] := \{j, \ldots, k\}$. For $n \in \Zset$, set
$\Zset_{\geq n} = \{n, n+1, \ldots\}$ and $\Rset_{\geq n} = [n, \infty)$.
The center of a Poisson algebra $(\PP, \{.,.\})$ will be denoted by 
\[
\ZZ(\PP) = \{ z \in \PP \mid \{z, f \} = 0, \forall f \in \PP \}.
\]
\medskip
\\
\noindent
{\bf Acknowledgements.} We would like to thank Ken Goodearl, Sergey Fomin, Jiang-Hua Lu, Bach Nguyen, Misha Shapiro 
and Kurt Trampel for helpful discussions and correspondence. M.Y. is grateful to Newcastle University and the Max Planck Institute 
for Mathematics in Bonn where part of this work was carried out.
\sectionnew{A rigidity Theorem}
\lb{Rig}
\subsection{Statement of main result} 
\label{2.1}
Let $\KK$ be a field of characteristic 0 and $T$ be an $m$-dimensional $\KK$-torus. Denote the
coordinate ring and rational function field of $T$ by
\[
\KK[T] \cong \KK[y_1^{\pm 1}, \ldots, y_m^{\pm 1}], \quad 
\KK(T) \cong \KK(y_1^{\pm 1}, \ldots, y_m^{\pm 1}). 
\]
For $\al=(i_1, \ldots, i_m) \in \Zset^m$, set 
\[
y^\al := y_1^{i_1} \ldots y_m^{i_m}. 
\]
Let 
\begin{equation}
\label{stand-basis}
\{\de_1, \ldots, \de_m\} \; \; \mbox{be the standard basis of $\Zset^m$}, 
\end{equation}
so $y^{\de_j} = y_j$. The algebra $\KK[T]$ is $\Zset^m$-graded by
\begin{equation}
\label{zm-grad}
\deg y^\al = \al, \quad \al \in \Zset^m.
\end{equation}

Fix a skewsymmetric additive bicharacter
\[
\Om \colon \Zset^m \times \Zset^m \to \KK
\]
and consider the quadratic Poisson bracket on $\KK[T]$ and $\KK(T)$ defined by 
\begin{equation}
\label{Om-P}
\{y^\al, y^\be \}_{\Om} := \Om(\al,\be) y^\al y^\be, \quad \forall \al,\be \in \Zset^m.
\end{equation}
It is graded of degree 0 with respect to the $\Zset^m$-grading \eqref{zm-grad}. 
In coordinates, this bracket is given by 
\[
\{y_k, y_j \}_{\Om} = \Om(\de_k,\de_j) y_k y_j, \quad \forall k,j \in [1,m].
\]
The Poisson structure $\{., . \}_{\Om}$ is called quadratic because it comes from the 
quadratic Poisson bivector field 
\[
\sum_{k,j} \Om(\delta_k, \delta_j) y_k y_j \partial_{y_k} \wedge \partial_{y_j}
\]
on $T$. The center of the Poisson algebra $(\KK[T], \{.,.\}_\Om)$ is
\begin{equation}
\label{cent}
\ZZ(\KK[T]) = \Span \{ y^\al \mid \al \in \rad \Om \},
\end{equation}
where
\[
\rad \Om = \{ \al \in \Zset^m \mid \Om(\al,\be) = 0, \forall \be \in \Zset^m \}.
\]
If $\be_1, \ldots \be_l$ is a free generating set of the (free) abelian group $\rad \Om$, then 
\[
\ZZ(\KK[T]) = \KK[ y^{\be_1}, \ldots, y^{\be_l}] \quad \mbox{and} \quad \ZZ(\KK(T)) = \KK(y^{\be_1}, \ldots, y^{\be_l}).
\]  

Let $D:=(d_1, \ldots, d_m) \in \Zset_{> 0}^m$ be a positive integral vector. The $\Zset^m$-grading \eqref{zm-grad} of $\KK[T]$ specializes
to a $\Zset$-grading given by 
\begin{equation}
\label{grad}
\deg y_k = d_k. 
\end{equation}
Its graded components will be denoted by $K[T]_n$, $n \in \Zset$. Set $K[T]_{\geq n} = \oplus_{l \geq n} K[T]_l$.
 
\bde{bifin-unip}
(i) An automorphism $\phi$ of the Poisson function field $(\KK(T), \{.,.\}_\Om)$ will be called {\em{bi-integral}} if
\[
\phi(y_k), \phi^{-1}(y_k) \in \KK[T], \quad \forall k \in [1,m].
\]

(ii) A bi-integral automorphism $\phi$ of $(\KK(T), \{.,.\}_\Om)$ will be called {\em{unipotent}} 
if 
\[
\phi(y_k) - y_k \in \KK[T]_{> d_k}, \quad \forall k \in [1,m].
\]
\ede

\bre{unip}
(i) It is easy to see that, if $\phi$ is a unipotent bi-integral automorphism of the Poisson function field $(\KK(T), \{.,.\}_\Om)$, 
then $\phi^{-1}$ has the same property; in particular, 
\[
\phi^{-1}(y_k) - y_k \in \KK[T]_{> d_k}, \quad \forall k \in [1,m].
\]

(ii) If $\phi$ and $\psi$ are unipotent bi-integral automorphisms of $(\KK(T), \{.,.\}_\Om)$, then the composition 
$\phi \circ \psi$ need not have the same property. That is, the set of unipotent bi-integral automorphisms of the Poisson function 
field $(\KK(T), \{.,.\}_\Om)$ does not form a subgroup of the group of all automorphisms of $(\KK(T), \{.,.\}_\Om)$. 
\ere
 
The next result is the main rigidity theorem in the paper.
 
\bth{1} Let $\KK$ be a field of characteristic 0, $m \in \Zset_{>0}$, $(d_1, \ldots, d_m) \in \Zset_{>0}^m$, 
and $\Om \colon \Zset^m \times \Zset^m \to \KK$ be a skewsymmetric additive bicharacter. Every 
unipotent bi-integral automorphism $\phi$ of the Poisson function field $(\KK(T), \{.,.\}_\Om)$, satisfies
\[
\phi(y_k) y_k^{-1} \in \ZZ(\KK[T]), \quad \forall k \in [1,m].
\]
\eth 
\subsection{An equivalent formulation}
\label{2.2}
\thref{1} has an equivalent formulation in terms of automorphisms of completions of Laurent polynomial rings.
The grading \eqref{grad} gives rise to the valuation of $\KK[T]$
\[
\nu \colon \KK[T] \to \Zset \sqcup \{ \infty\}, \quad
\nu(r_j + \cdots + r_k) = j, \forall r_i \in \KK[T]_i, i \in [j,k], r_j \neq 0.
\]
The corresponding completion of $\KK[T]$ is
\[
\wh{\KK[T]} := \{ u_n + u_{n+1} + \cdots \mid n \in \Zset,  u_l \in \KK[T]_l \}.
\] 
It has the descending $\Zset$-filtration, defined by
\[
\wh{\KK[T]}_{\geq n} := \{ u_n + u_{n+1} + \cdots \mid  u_l \in \KK[T]_l \} \quad \mbox{for} \; \; n \in \Nset.
\]
The Poisson structure \eqref{Om-P} has a canonical extension to  $\wh{\KK[T]}$ which will be also denoted by $\{.,.\}_\Om$. 

\bth{2} Let $\KK$ be a field of characteristic 0, $m \in \Zset_{>0}$, $(d_1, \ldots, d_m) \in \Zset_{>0}^m$, 
and $\Om \colon \Zset^m \times \Zset^m \to \KK$ be a skewsymmetric additive bicharacter. For every continuous automorphism $\theta$ 
of the Poisson algebra $(\wh{\KK(T)}, \{.,.\}_\Om)$, satisfying
\begin{enumerate}
\item[(*)] $\theta(y_k) - y_k, \theta^{-1}(y_k) - y_k \in \KK[T]_{> d_k}$, $\forall k \in [1,m]$,
\end{enumerate}
we have
\[
\theta(y_k) y_k^{-1} \in \ZZ(\KK[T]), \quad \forall k \in [1,m].
\]
\eth
Theorems  \ref{t1} and \ref{t2} are equivalent due to the following lemma. 
\ble{equiv} There is a bijection between the set of unipotent bi-integral automorphisms $\phi$ of the Poisson function field $(\KK(T), \{.,.\}_\Om)$ 
and the set of continuous automorphisms $\theta$ of $(\wh{\KK[T]}, \{.,.\}_\Om)$ satisfying the condition {\em{(*)}} in \thref{2}. 
The bijection is uniquely defined by the condition that 
\[
\phi|_{\KK[T]} = \theta|_{\KK[T]}.
\]
\ele
\begin{proof} Let $\phi$ be a unipotent bi-integral automorphism of the Poisson function field $(\KK(T), \{.,.\}_\Om)$. Thus,
\[
\phi(y_k) = (1+f_k) y_k \; \; \mbox{for some} \; \; f_k \in \KK[T]_{\geq 1}.
\]
Define 
\[
\theta(y_k) := \phi(y_k)=  (1+f_k) y_k, \; \; \theta(y_k^{-1} ) = (1+ f_k)^{-1} y_k^{-1} := \sum_{n=0}^\infty (-1)^n f_k^n y_k^{-1}  \in \wh{\KK[T]}.
\]
It is easy to see that $\theta$ uniquely extends to a continuous automorphism of $(\wh{\KK[T]}, \{.,.\}_\Om)$ 
which satisfies the condition (*) in \thref{2}. 

In the opposite direction, let $\theta$ be a continuous automorphism of the Poisson algebra $(\wh{\KK[T]}, \{.,.\}_\Om)$ which satisfies the condition (*) in \thref{2}. 
Define
\[
\phi(y_k) := \theta(y_k), \; \; \phi(y_k^{-1}) := \theta(y_k)^{-1} \in \KK(T). 
\]
It is clear that $\phi$ uniquely extends to a unipotent bi-integral automorphism of the Poisson function field $(\KK(T), \{.,.\}_\Om)$.
\end{proof}
\subsection{Support and cone of unipotent bi-integral automorphisms of the Poisson function field $(\KK[T], \{.,.\}_\Om)$}
\label{2.3}
By a {\em{strict polyhedral cone}} in $\Rset^m$ we will mean a cone of the form
\[
C= \Rset_{\geq 0} X = \{ r_1 \al_1 + \cdots + r_n \al_n \mid r_i \in \Rset_{\geq 0}, \al_i \in X \}
\]
for a finite subset $X \subset \Rset^m$ such that $\al \in C \Rightarrow - \al \notin C$ for all $\al \in \Rset^m$, $\al \neq 0$. A ray $R_{\geq 0} \al$ 
in $\Rset^m$ will be called {\em{extremal}} for $C$ if for all $\al_1, \al_2 \in C$,
\[
\al_1 + \al_2 \in \Rset_{\geq 0} \al \Rightarrow \al_1, \al_2 \in \Rset_{\geq 0} \al.
\]

Given 
\[
f = \sum_{\be \in \Zset^m} c_\be y^\be \in \KK[T], c_\be \in \KK
\]
and $\al \in \Zset^m$, set
\begin{equation}
\label{coeff}
[f]_\al := c_\al.
\end{equation}
\bde{supp}
(i) Define the {\em{support}} of an element $f \in \KK[T]$ to be the set 
\[
\Supp(f) = \{ \al \in \Zset^m \mid [f]_\al \neq 0 \}.
\]
(ii) For a unipotent bi-integral automorphism $\phi$ of the Poisson function field $(\KK[T], \{.,.\}_\Om)$, the set
\[
\Supp(\phi) := \Big( \Nset \, \cdot \,  \bigcup_j \Supp(\phi(y_j)y_j^{-1} -1) \Big) \backslash \{ 0\}
\]
will be called the {\em{support}} of $\phi$. 
\ede
The following fact was proved in \cite{Y-AD} in the greater generality of automorphisms of completions of quantum tori.
\ble{Supp} \cite[Lemma 3.10]{Y-AD} Let $\phi$ be a unipotent bi-integral automorphism of $(\KK(T), \{.,.\}_\Om)$. Then 
\[
\Supp( \phi^{-1} ) = \Supp (\phi), 
\]
and for all $\al \in \Zset^m$, 
\begin{equation}
\label{supp-al}
\Supp( \phi(y^\al)  - y^\al ) \subset \al + \Supp(\phi).
\end{equation}

Furthermore, if $\psi$ is another unipotent bi-integral automorphism of $(\KK(T), \{.,.\}_\Om)$ such that $\phi \circ \psi$ 
has the same property, then
\[
\Supp (\phi \circ \psi) \subseteq \Big( (\Supp(\phi) \cup \{0\}) + ( \Supp(\psi) \cup \{0\} ) \Big) \backslash \{ 0 \},
\]
recall \reref{unip} (ii).
\ele

Denote the functional
\begin{equation}
\label{muD}
\mu_D \colon \Rset^m \to \Rset, \quad \mu(r_1, \ldots, r_m) = r_1 d_1 + \cdots + r_m d_m, \forall k_i \in \Rset.
\end{equation}
\bde{cone}
Define the {\em{cone}} of a unipotent bi-integral  automorphism $\phi$ of $(\KK(T), \{.,.\}_\Om)$ to be the set
\[
\Con(\phi) = \Rset_{\geq 0} \Supp (\phi) =  \{ r_1 \al_1 + \cdots + r_n \al_n \mid r_i \in \Rset_{\geq 0}, \al_i \in \Supp(\phi) \}.
\]
\ede
Since the support of each unipotent bi-integral automorphism of $(\KK(T), \{.,.\}_\Om)$ satisfies
\[
\Supp(\phi) \subset \{ \al \in \Zset^m \mid \mu_D(\al) \geq 1 \}
\]
and $d_1, \ldots, d_m \in \Zset_{>0}$, we obtain the following:
\bco{con}
The cone of each unipotent bi-integral automorphism $\phi$ of the Poisson function field $(\KK(T), \{.,.\}_\Om)$ is a strict polyhedral cone. 
Furthermore, the set $\Con(\phi) \backslash \{ 0 \}$ is contained in the strict half-space
\[
\{ \al \in \Rset^m \mid \mu_D(\al) >0  \}.
\]
\eco
\subsection{Restrictions of unipotent bi-integral automorphisms to extremal rays}
\label{2.4}
The following proposition will play a key role in the proof of \thref{1}. For a subset $X \subset \Rset^m$ and $f \in K[T]$, denote 
\[
f|_X := \sum_{\be \in X \cap \Zset^m}  [f]_\be y^\be,
\]
using the notation in \eqref{coeff}.
\bpr{restr} Assume that $\phi$ is a unipotent bi-integral automorphism of the Poisson function field $(\KK(T), \{.,.\}_\Om)$ and that $\Rset_{\geq 0} \al$ is an 
extremal ray of $\Con(\phi)$. Then 
\begin{equation}
\label{def-restr}
\phi|_{\Rset_{\geq 0} \al} (y^\be) := \phi(y^\be)|_{\be + \Rset_{\geq 0} \al}  \quad \mbox{for} \; \; \be \in \Zset_{>0}^m
\end{equation}
uniquely defines a unipotent bi-integral automorphism of $(\KK(T), \{.,.\}_\Om)$. Its inverse is $(\phi^{-1})|_{\Rset_{\geq 0} \al}$.
\epr
The automorphism of $(\KK(T), \{.,.\}_\Om)$ will be called {\em{the restriction of $\phi$ to the extremal ray $\Rset_{\geq 0} \al$}}. By its definition,
\[
\Supp \big( \phi|_{\Rset_{\geq 0} \al} \big) = \Supp(\phi) \cap \Rset_{\geq 0} \al.
\]
\begin{proof} It follows from Proposition 3.11 (i) in \cite{Y-AD} that \eqref{def-restr} defines an automorphism of $\KK(T)$. The second
part of Proposition 3.11 in \cite{Y-AD} implies that the inverse of this automorphism is $(\phi^{-1})|_{\Rset_{\geq 0} \al}$. These facts can 
be also deduced from \leref{Supp}, we leave the details to the reader. The fact that $\phi|_{\Rset_{\geq 0} \al}$ is unipotent and bi-integral, and the uniqueness 
statement in \prref{restr} both follow directly. It remains to prove that $\phi|_{\Rset_{\geq 0} \al}$ is an automorphism of Poisson function fields, for which it is sufficient 
to prove that
\[
\phi |_{\Rset_{\geq 0} \al} : \Span \{ y^\al \mid \al \in \Zset_{>0}^m\}  \to \KK[T]
\]
is a homomorphism of Poisson algebras with respect to the restrictions of the Poisson structure $\{.,.\}_\Om$. Using the fact that the 
Poisson bracket $\{.,.\}_\Om$ has degree 0 in the $\Zset^m$-grading \eqref{zm-grad}, we obtain that for all $\be, \ga \in \Zset_{>0}^m$, 
\begin{align*}
&\{ \phi|_{\Rset_{\geq 0} \al} (y^\be), \phi|_{\Rset_{\geq 0} \al} (y^\ga) \}_\Om = 
\{ \phi(y^\be)|_{\be+ \Rset_{\geq 0} \al} , \phi(y^\ga)|_{\ga + \Rset_{\geq 0} \al}\}_\Om \\
=& \{ \phi(y^\be), \phi(y^\ga)\}_\Om|_{\be+ \ga + \Rset_{\geq 0} \al} = \phi( \{ y^\be, y^\ga\}_\Om) |_{\be+ \ga + \Rset_{\geq 0} \al} = 
\phi|_{\Rset_{\geq 0} \al} ( \{ y^\be, y^\ga \}_\Om). 
\end{align*}
In the second equality we also used the assumption that $\Rset_{\geq 0} \al$ is an extremal ray of $\Supp(\phi)$ and eq. \eqref{supp-al} in \leref{Supp}. 
\end{proof} 
Let $\phi$ be a unipotent bi-integral automorphism of the Poisson function field $(\KK(T), \{.,.\}_\Om)$ and $\theta$ be the continuous automorphism of 
the Poisson algebra $(\wh{\KK[T]}, \{.,.\}_\Om)$ corresponding to it under the bijection from \leref{equiv}. Then the continuous automorphism of $(\wh{\KK[T]}, \{.,.\}_\Om)$, corresponding to the restriction $\phi|_{\Rset_{\geq 0} \al}$ is 
given by
\[
\theta|_{\Rset_{\geq 0} \al} (y^\be) := \theta(y^\be)|_{\be + \Rset_{\geq 0} \al}  \quad \mbox{for} \; \; \be \in \Zset^m.
\]
\subsection{Proof of \thref{1}}
\label{2.5}
We will call a continuous automorphism $\theta$ of $(\wh{\KK[T]}, \{.,.\}_\Om)$ {\em{unipotent}} if
\[
\theta(y_j) - y_j \in \wh{\KK[T]}_{> d_j}, \quad \forall j \in [1,m].
\]
It is easy to see that this is equivalent to saying that 
\[
\theta(y^\be) - y^\be \in \wh{\KK[T]}_{> \mu_D(\be)}, \quad \forall \be \in \Zset^m
\]
in terms of the functional in \eqref{muD}.

\ble{der} Let $\al \in \Zset$, $\al \notin \rad \Om$.

(i) If $\partial$ is a $\Zset^m$-graded derivation of the Poisson algebra $(\KK[T], \{.,.\}_\Om)$
of degree $\al$ {\em{(}}with respect to the grading \eqref{zm-grad}{\em{)}}, then 
\[
\partial = a \{ y^\al, - \}_\Om, 
\]
for some $a \in \KK$.

(ii) If $\theta$ is a unipotent continuous automorphism of $(\wh{\KK[T]}, \{.,.\}_\Om)$ with support contained in $\Zset_{>0} \al$, 
then 
\begin{equation}
\label{Poisson-prod}
\theta = \prod_{n=1}^\infty  \exp ( a_n \{ y^{n \al}, - \}_\Om)
\end{equation}
for some $a_1, a_2, \ldots \in \KK$, where the product is taken in decreasing order from left to right.
\ele
Note that, since the Poisson structure $\{.,.\}_\Om$ is graded of degree 0 with respect to the $\Zset^m$-grading \eqref{zm-grad}, the 
composition in the RHS of \eqref{Poisson-prod}
\[
\ldots \exp ( a_2 \{ y^{2 \al}, - \}_\Om) \exp ( a_1 \{ y^{\al}, - \}_\Om)
\]
is a well defined unipotent continuous automorphism of $(\wh{\KK[T]}, \{.,.\}_\Om)$ for all $a_n \in \KK$.
\begin{proof} (i) Let $k \in [1,m]$ be such that $\Om(\al, \de_k) \neq 0$. Denote
\[
\partial(y_j) = b_j y^{\al+\de_j} \; \; \mbox{for some} \; \; b_1, \ldots, b_m \in \KK. 
\]
The derivation properties of $\partial$ with respect to the commutative product and the Poisson bracket, 
and the definition of $\{.,.\}_\Om$ give
\[
\partial(\{y_k, y_j\}_\Om) = \{\partial(y_k), y_j \}_\Om + \{y_k, \partial(y_j)\}_\Om = \left( b_k \Om(\al + \de_k, \de_j) + b_j \Om(\de_k, \al + \de_j) \right) y^{\al+\de_k+\de_j}
\]
and
\[
\partial(\{y_k, y_j\}_\Om) = \Om(\de_k, \de_j) \partial(y_k y_j) = (b_k +b_j) \Om(\de_k, \de_j) e^{\al+\de_k+\de_j}. 
\]
This implies that 
\[
b_j = \frac{\Om(\al,\de_j)}{\Om(\al,\de_k)} b_k
\]
for all $j \in [1,m]$, and thus,
\[
\partial = \frac{b_k}{\Om(\al,\de_k)} \{y^\al, -\}_\Om.
\]

(ii) For $j \in [1,m]$, define $\partial(y_j) \in \KK[T]_{\al + \de_j}$ by 
\[
\partial(y_j) := [\theta(y_j)]_{\al+ \de_j} y^{\al+ \de_j}.
\] 
Since 
\[
\Supp( \phi(y_j) - y_j - \de(y_j)) \in \de_j + \Zset_{\geq 2} \al,
\]
and $\phi$ an automorphism of the Poisson algebra $(\wh{\KK[T]}, \{.,.\}_\Om)$, 
$\partial$ extends to a derivation of $(\KK[T], \{.,.\}_\Om)$. 
By part (i), $\partial = a_1 \{ y^\al, - \}_\Om$ for some $a_1 \in \KK$. Hence, 
\[
\psi :=  \exp ( - a_1 \{ y^{\al}, - \}_\Om) \phi
\]
is a unipotent continuous automorphism of $(\wh{\KK[T]}, \{.,.\}_\Om)$ whose support is contained in $\Zset_{\geq 2} \al$. 
In the same way one shows that there exists $a_2 \in \KK$ such that 
\[
\exp ( a_2 \{ y^{2 \al}, - \}_\Om) \psi 
\]
is a unipotent continuous automorphism of $(\wh{\KK[T]}, \{.,.\}_\Om)$, whose support is contained in $\Zset_{\geq 2} \al$. 
The proof of part (ii) is now completed by induction.
\end{proof}
\bco{fix} Let $\al \in \Zset$, $\al \notin \rad \Om$. For every unipotent bi-integral automorphism $\phi$ of the Poisson function field 
$(\KK(T), \{.,.\}_\Om)$ with $\Supp(\phi) \subseteq \Zset_{\geq 1} \al$, we have 
\[
\phi(y^{\al}) = y^\al.
\]
\eco
The corollary follows from \leref{der} (ii) and the fact that for every unipotent bi-integral automorphism of 
$(\KK(T), \{.,.\}_\Om)$, the continuous automorphism of $(\wh{\KK[T]}, \{.,.\}_\Om)$, corresponding to it under \leref{equiv}, is unipotent.
\medskip
\\
\noindent
{\bf{Proof of \thref{1}.}} Because of \eqref{cent}, the statement is equivalent to the inclusion
\[
\Supp(\phi) \subset \rad \Om.
\]
Assume that this is not the case. Since $\rad \Om$ is the intersection of a subspace of $\Rset^m$ with $\Zset^m$, 
$\Con(\phi)$ will have an extremal ray $\Rset_{\geq 0} \al$ such that 
\begin{equation}
\label{alnot}
\al \notin \rad \Om.
\end{equation}
By rescaling $\al$ we can assume that 
\[
\Rset_{\geq 0} \al \cap \Zset^m = \Nset \al.
\]
Denote
\[
\psi:= \phi|_{\Rset_{\geq 0} \al}.
\]
By the definition of support of $\phi$, there exists $k \in [1,m]$ such that 
\[
\psi(y_k) = \sum_{i=0}^n a_i y^{i \al} y_k \quad \mbox{with} \; \; a_i \in \KK
\]
for some $n>0$, $a_n \neq 0$. It follows from \prref{restr} that $\psi$ is a unipotent 
bi-integral automorphism of $(\KK(T), \{.,.\}_\Om)$ and $\psi^{-1} = \phi^{-1}|_{\Rset_{\geq 0 } \al}$. 
So, 
\[
\psi^{-1}(y_k) = \sum_{j=0}^l a'_j y^{j \al} y_k \quad \mbox{with} \; \; a'_j \in \KK
\]
for some $l \geq 0$, $a'_l \neq 0$. By applying \coref{fix}, we obtain 
\[
y_k = \psi \psi^{-1} (y_k) =  \sum_{i=0}^n \sum_{j=0}^l a_i a'_j y^{(i+j) \al} y_k. 
\]
However, $n+l >0$ and the coefficient of $y^{(n+l)\al + \de_k} = y^{(n+l) \al} y_k$ in the LHS of the last equality is 0, 
while in the coefficient of $y^{(n+l) \al} y_k$ in the RHS is $a_n a'_l \neq 0$. This is a contradiction, which completes the proof of the theorem.
\qed
\sectionnew{Poisson automorphisms of cluster algebras}
\lb{App-cl}
\subsection{Poisson clusters and cluster algebras}
\bde{One-cl} Assume that $(\PP, \{.,.\})$ is a Poisson algebra over a field $\KK$. We will say that the $n$-tuple $(y_1, \ldots, y_m)$ 
of elements of  $\PP$ is a 
{\em{Poisson cluster}} of $\PP$, if the following three conditions are satisfied:
\begin{enumerate}
\item $\{y_j, y_k \} = c_{jk} y_j y_k$, $\forall j,k \in [1,m]$ for some $c_{jk} \in \KK$. 
\item $y_1, \ldots, y_n$ generate a polynomial subalgebra of $\PP$. 
\item $\PP$ is contained in the localization $\KK[y_1^{\pm1}, \ldots, y_m^{\pm 1}]$.
\end{enumerate} 
\ede
In particular, each Poisson cluster of $(\PP, \{.,.\})$ gives rise to the embeddings of Poisson algebras
\begin{equation}
\label{embed}
\KK[y_1, \ldots, y_m] \subseteq \PP \subseteq \KK[y_1^{\pm1}, \ldots, y_m^{\pm 1}].
\end{equation} 
The first algebra is a Poisson subalgebra of $(\PP, \{.,.\})$ because of the condition (1). Every localization of a Poisson algebra has a canonical structure of
Poisson algebra. Hence, the Laurent polynomial ring $\KK[y_1^{\pm1}, \ldots, y_m^{\pm 1}]$ admits a canonical Poisson algebra structure as a localization of the first 
algebra. Because of the embeddings in \eqref{embed}, the third algebra is also a localization of the second, namely, 
\[
\KK[y_1^{\pm1}, \ldots, y_m^{\pm 1}] \cong \PP[y_1^{-1}, \ldots, y_m^{-1}]
\]
as Poisson algebras.

Seeds of cluster algebras give rise to Poisson clusters as follows. Assume that $\AA$ is {\em{a cluster algebra of geometric type}} \cite{FZ1} with 
base ring extended from $\Zset$ to the field $\KK$. Let $\Sig = (y_1, \ldots, y_m, \wt{B})$ be a {\em{labeled seed}} of $\AA$, 
i.e., one whose cluster variables $y_1, \ldots, y_m$ are labelled with the integers in $[1,m]$. Here $\wt{B}$ is {\em{the exchange matrix}} of $\Sig$; it is
an integral matrix of size $m \times n$ whose principal $n \times n$ submatrix $B$ is skew-symmetrizable (i.e., 
$\diag(d_1, \ldots d_n) B$ is a skewsymmetric matrix for some collection $d_1, \ldots, d_n$ of positive, relative prime integers).   
The number $n$ is the number of {\em{exchangeable cluster variables}} of $\AA$. We will assume throughout 
that $y_{n+1}, \ldots, y_m$ are {\em{the frozen variables}} of $\AA$. We refer the reader to \cite{FZ1,GSVb,W} for details on cluster algebras.

\bde{comp} \cite{BZ,GSVb} Let $\La = (\la_{jk})$ be a skew symmetric integral matrix. The pair $(\La, \wt{B})$ is called compatible if 
the $n\times n$ principle submatrix of the $m \times n$ matrix $- \La \wt{B} = \La^t \wt{B} $ equals $\diag(d_1, \ldots, d_n)$ and all other entries of it vanish.
\ede
Here and below $X^t$ denotes the transpose of a matrix $X$.
Following \cite{BZ,GSV0}, for a mutable index $k \in [1,n]$, define the $m \times m$ matrix
\[
E_k=(e_{ij}), \; \; e_{ij} = 
\begin{cases}
\de_{ij}, & \mbox{if} \; \: j \neq k \\
-1, & \mbox{if} \; \; i= j =k \\
\max(0, b_{ik}), & \mbox{if} \; \; i\neq j = k
\end{cases}
\]
where $b_{ij}$, $i \in [1,m]$, $j \in [1,n]$ are the entries of the exchange matrix $\wt{B}$. 
If $(\La, \wt{B})$ is a compatible pair and $k \in [1,n]$, then \cite{BZ,GSV0} the matrix 
\[
\La' = E^t_k \La E_k
\]  
has the property that for the seed mutation $\mu_k(\Sig) = (y_1, \ldots, y_{k-1}, y'_k, y_{k+1}, \ldots, y_n, \mu_k(\wt{B}))$, the pair 
$(\La', \mu_k({\wt{B}}))$ is compatible. (We refer the reader to \cite[Sect. 3]{BZ} for a detailed study on the properties of compatible pairs.)
Iterating this formula, simultaneously with the seed mutation formulas, 
one constructs compatible skewsymmetric integral matrices with all seeds 
of the algebra $\AA$, and defines the Gekhtman--Shapiro--Vainshtein Poisson structure of $\AA$ (associated to $\La$) to be given 
by 
\begin{equation}
\label{GSV}
\{y_j, y_k\} = \la_{jk} y_j y_k, \quad \forall j,k \in [1,m]. 
\end{equation}
More precisely, by the Laurent phenomenon \cite{FZ-L}, $\AA \subseteq \KK[y_1^{\pm 1}, \ldots, y_m^{\pm 1}]$. It was proved in \cite{GSV0}
that $\AA$ is a Poisson subalgebra of  $\KK[y_1^{\pm 1}, \ldots, y_m^{\pm 1}]$, equipped with the Poisson bracket \eqref{GSV}, and that 
the Poisson bracket between the cluster variables in every seed of $\AA$ is given by the analogous to \eqref{GSV} in terms of the 
associated $\La$-matrix to the (labeled) seed.  This formula implies the following:
\bco{sl} Let $\AA$ be a cluster algebra of geometric type equipped with a Gekhtman--Shapiro--Vainshtein Poisson structure $\{.,.\}$. 
For every labeled seed $\Sig = (y_1, \ldots, y_m, \wt{B})$ of $\AA$, $(y_1, \ldots, y_n)$ is a Poisson cluster of $(\AA, \{.,.\})$.
\eco
\subsection{Automorphisms of $\Nset$-graded connected algebras admitting a homogeneous cluster}
\label{3.2}
Let $(\PP= \oplus_{n \in \Nset} \PP_n , \{.,.\})$ be an {\em{$\Nset$-graded Poisson algebra}}. The Poisson bracket will be assumed to be graded 
but not necessarily of degree 0. Denote by $\Aut(\PP)$ the group of automorphisms of the Poisson algebra $(\PP, \{.,.\})$  and by
\[
\UAut (\PP) = \{ \phi \in \Aut (\PP) \mid \phi(f) - f \in \PP_{\geq n+1}, \; \; \forall n \in \Nset, f \in \PP_n \}
\]
the subgroup of {\em{unipotent automorphisms}} of $(\PP, \{.,.\})$. When $\dim \PP_0$ is small (especially when $\PP$ is connected, i.e., $\dim \PP_0 =1$), 
it is easy to describe $\Aut(\PP)$ in terms of $\UAut(\PP)$. Thus, the 
problem for determining $\Aut(\PP)$ reduces to that of determining $\UAut(\PP)$. 
The latter problem is extremely hard and has been solved only in very few situations.

Our next result proves that all unipotent automorphisms of $\Nset$-graded 
connected algebras admitting a homogeneous cluster, 
have a very special form that depends only on the Poisson center of the associated Laurent polynomial ring. It can be 
used to fully describe the automorphism groups of such algebras; in the next section this is illustrated with the coordinate ring of the open Schubert cells of all full flag 
varieties equipped with the standard Poisson structure. 
\bth{3} Assume that $(\PP= \oplus_{n \in \Nset} \PP_n , \{.,.\})$ is an $\Nset$-graded Poisson algebra over a field 
$\KK$ of characteristic 0 and $(y_1, \ldots, y_m)$ is a Poisson cluster of $\PP$ 
consisting of homogeneous elements of positive degrees. Then every unipotent Poisson automorphism $\phi$ of $(\PP, \{.,.\})$ has the property 
\[
\phi(y_k) y_k^{-1} \in \ZZ( \KK[y_1^{\pm 1}, \ldots, y_m^{\pm 1}]), \quad \forall k \in [1,m],
\]
where $\ZZ(.)$ refers to the center of the Poisson Laurent polynomial ring with the Poisson structure $\{.,.\}$. 
\eth
\thref{3} places a very strong restriction on the possible form of the unipotent Poisson automorphisms of $(\PP, \{.,.\})$ because 
in all important situations the Krull dimension of $\ZZ( \KK[y_1^{\pm 1}, \ldots, y_m^{\pm 1}])$ is much smaller that that 
of $\PP$, see next section. Furthermore, because of the second embedding in \eqref{embed}, 
$\phi$ is fully determined from the values $\phi(y_k)$ for $k \in [1,m]$.

Assume the setting of \thref{3} and denote by $\Om \colon \Zset^m \times \Zset^m \to \KK$ the skewsymmetric additive bicharacter given by 
$\Om(\de_j, \de_k) := \la_{jk}$. We have the isomorphism of Poisson algebras
\begin{equation}
\label{is}
(\KK[y_1^{\pm 1}, \ldots, y_m^{\pm 1}], \{.,.\}) \cong (\KK[T], \{.,.\}_\Om)
\end{equation}
where the first algebra is equipped with the restricted Poisson structure from $\PP$ and the second is the 
one considered in Sect. \ref{Rig}. Denote by $d_k \in \Zset_{>0}$ the degree of $y_k$.  

\thref{3} directly follows from \thref{1} and the following proposition.

\bpr{aut-emb} In the setting of \thref{3}, we have a group embedding 
\[
\iota \colon \UAut (\PP, \{.,.\}) \hra \Aut(\KK[T], \{.,.\}_\Om)
\]
whose image is contained in the set of unipotent bi-integral automorphisms of $(\KK[T], \{.,.\}_\Om)$.
For $\phi \in \UAut (\PP, \{.,.\})$, $\iota(\phi)$ is the unique automorphism of $(\KK[T], \{.,.\}_\Om)$ such that
\[
\iota(\phi)|_\PP = \phi. 
\]
\epr
\begin{proof} Denote $\KK[\Aset^m]:= \KK[y_1, \ldots, y_m]$. It is a Poisson 
subalgebra of $(\KK[T], \{,.\}_\Om)$. Identifying the Poisson algebras in \eqref{is} and taking into account \eqref{embed}, gives the embeddings 
of Poisson algebras
\begin{equation}
(\KK[\Aset^m], \{,.\}_\Om) \subseteq (\PP, \{,.\}) \subseteq (\KK[T], \{,.\}_\Om).
\label{embed2}
\end{equation}
Since the elements $y_1, \ldots, y_m$ are homogeneous, these embeddings are graded. Fix $\phi \in \UAut (\PP, \{.,.\})$ and define 
\[
\Phi \colon \KK[\Aset^m] \to \KK[T] \quad \mbox{by} \; \; \Phi(f) := \phi(f) \in \PP \subseteq \KK[T], \; \forall 
f \in \KK[\Aset^m] \subseteq \PP.
\]
Taking into account that $\phi$ is an automorphism of $(\PP, \{,.\})$ and using the embeddings in \eqref{embed2}, we obtain that
$\Phi$ a homomorphism of Poisson algebras with respect to the Poisson structure $\{.,.\}_\Om$. Since 
\begin{equation}
\label{Phi-unip}
\Phi(f) - f \in \PP_{\geq n+1} \subseteq \KK[T]_{\geq n +1}, \forall f \in \KK[\Aset^m]_n \subseteq \PP_n, n \in \Nset,
\end{equation}
$\Phi$ is injective. Thus, it extends to an injective Poisson homomorphism  $\Phi \colon (\KK(T), \{.,.\}_\Om) \to (\KK(T), \{.,.\}_\Om)$. 
Since $\PP$ is a subalgebra of $\KK[T]$, it is an integral domain. This, the definition of $\Phi$ and the fact that $\PP$ is 
contained in the field of fractions of $\KK[\Aset^m]$ imply that 
\begin{equation}
\label{phiPhi}
\Phi|_\PP = \phi.
\end{equation}
Denote by $\Psi \colon (\KK(T), \{.,.\}_\Om) \to (\KK(T), \{.,.\}_\Om)$ the injective Poisson homomorphism obtained in the same way 
from $\phi^{-1} \in \UAut (\PP, \{.,.\})$. It follows from \eqref{phiPhi} that
\[
\Psi \Phi(f) = \Psi \phi(f) = \phi^{-1} \phi (f) = f, \quad \forall f \in \KK[\Aset^m]
\]
because $\phi(f) \in \PP$. Therefore, $\Phi$ is a Poisson automorphism of $(\KK(T), \{.,.\}_\Om)$ with inverse $\Psi$, and by \eqref{Phi-unip}, 
$\Phi$ is a unipotent bi-integral automorphism. Define $\iota(\phi) := \Phi$. Eq. \eqref{phiPhi} implies that 
$\iota \colon \UAut (\PP, \{.,.\}) \to\Aut(\KK[T], \{.,.\}_\Om)$ is injective. We have proved that its image is contained in 
unipotent bi-integral automorphisms of $(\KK[T], \{.,.\}_\Om)$. It remains to prove that $\iota$ is a group homomorphism.

Let $\phi_1, \phi_2 \in \UAut (\PP, \{.,.\})$. Applying \eqref{phiPhi} gives 
\[
\big( \iota(\phi_1) \circ \iota(\phi_2) \big) (f) = \iota(\phi_1) \big( \phi_2 (f) \big) = \big( \phi_1 \circ \phi_2 \big) (f) =
\iota(  \phi_1 \circ \phi_2) (f), \quad \forall f \in \KK[\Aset^m]
\]
because $\phi(f) \in \PP$. Since $\PP$ is contained in the field of fractions of $\KK[\Aset^m]$, we obtain that 
$\iota(\phi_1) \circ \iota(\phi_2) = \iota(  \phi_1 \circ \phi_2)$, which completes the proof of the proposition
\end{proof}
\subsection{Poisson automorphisms of $\Nset$-graded connected cluster algebras} 
\label{3.3}
\thref{3} and \coref{sl} imply the following result. 
\bth{4} Let $\AA$ be an $\Nset$-graded cluster algebra of geometric type over a field $\KK$ of characteristic 0, 
equipped with a Gekhtman--Shapiro--Vainshtein Poisson structure $\{.,.\}$. Let $\Sig = (y_1, \ldots, y_m, \wt{B})$ be a labeled 
seed of $\AA$ such that $y_1, \ldots, y_n$ are homogeneous of positive degrees. For every unipotent automorphism $\phi$ 
of the Poisson algebra $(\AA, \{.,.\})$, 
\[
\phi(y_k) y_k^{-1} \in \ZZ( \KK[y_1^{\pm 1}, \ldots, y_m^{\pm 1}])
\]
where the center is computed with respect to the induced Poisson structure on the localization 
$\KK[y_1^{\pm 1}, \ldots, y_m^{\pm 1}] = \AA[y_1^{-1}, \ldots, y_m^{-1}]$.
\eth

There are two important subgroups of the automorphism group of the Poisson cluster algebra $(\AA, \{.,.\})$: the subgroup of {\em{the modular Poisson automorphisms}} 
and the subgroup of {\em{toric transformations}}. 
For a labeled seed $\Sig$ of $\AA$, denote by $y_1(\Sig), \ldots, y_m(\Sig)$ 
the cluster variables in it, by $\wt{B}(\Sig)$ its exchange matrix, and by $\La(\Sig)$ the compatible integral skewsymmetric $m \times m$-matrix.
As before, we will assume that $y_{n+1}(\Sig), \ldots, y_m(\Sig)$ are the frozen variables of $\AA$.

{\bf{(1)}} The {\em{group of toric transformations}} of $(\AA, \{.,.\})$ was introduced by Gekhtman--Shapiro--Vainshtein in \cite{GSV0} as follows. 
For a labeled seed $\Sig$ of $\AA$, denote the lattice 
\[
\Theta_\Sig := \Ker \wt{B}(\Sig)^t \cap \Zset^m
\]
where $\wt{B}(\Sig)^t$ is the transpose of the matrix $\wt{B}(\Sig)$.
In \cite[Proposition 3.3]{BZ} it was proved that the assumption that $( \La(\Sig), \wt{B}(\Sig))$ is a compatible pair implies that
$\rank \wt{B}(\Sig) = n$, so $\rank \Theta_\Sig = m-n$. Denote by $T_\AA$, the $m-n$-dimensional 
$\KK$-subtorus of $(\KK^*)^m$ generated by 
\[
(a^{l_1}, \ldots, a^{l_m}) \quad \mbox{for} \; \;  (l_1, \ldots, l_m) \in \Theta_\Sig, a \in \KK^*.
\]
It was proved in \cite[Lemma 2.2]{GSV0} that the action of $T_\AA$ on $\KK[y_1(\Sig)^{\pm 1}, \ldots, y_m(\Sig)^{\pm 1}]$, given by
\[
(a_1, \ldots, a_m) \cdot y_k(\Sig) := a_k y_k(\Sig), \quad \forall k \in [1,m], (a_1, \ldots, a_m) \in T_\AA
\]
preserves $\AA$ and that, for any other labeled seed of $\AA$, the two tori and the corresponding actions are canonically isomorphic. 
The action of $T_\AA$ obviously preserves $\{.,.\}$ for homogeneity reasons, so, we have an embedding
\[
T_ \AA \subset \Aut(\AA, \{.,.\}).
\]
The group $T_\AA$ is called \cite{GSV0} {\em{group of toric transformations}} of $(\AA, \{.,.\})$.

{\bf{(2)}} The {\em{modular group}} of a cluster algebra $\AA$ of geometric type was defined in \cite{FZ2,FG,ASS} (a more general construction
for quasiautomorphsims in the presence of coefficients in semifields was introduced in \cite{Fr}). We use the term {\em{modular Poisson group}} 
for the subgroup of the former that preserves the Poisson structure $\{.,.\}$.

A {\em{modular Poisson automorphism}} of $\AA$ is defined to be any automorphism $\phi \in \Aut(\AA)$ 
that has the property that there is a cluster $(y_1, \ldots, y_m)$ of $\AA$ satisfying the conditions:
\begin{enumerate}
\item {\em{$(\phi(y_1), \ldots, \phi(y_m))$ is a cluster of $\AA$ and $\phi(y_j) = y_j$ for $j \in [n+1,m]$,}}
\item {\em{$\{\phi(y_k), \phi(y_j) \} = \phi(\{y_k,y_j\})$ for all $k,j \in [1,m]$,}}
\item {\em{$\phi$ commutes with mutations in the sense that $\phi(\mu_k (y_k)) = \mu_k (\phi(y_k))$ for all $k \in [1,n]$.}} 
\end{enumerate}
It is easy to check that this implies that every cluster of $\AA$ will possess the same properties.

Recall that an $m \times n$ matrix $\wt{B}$ is called indecomposable if there is no nontrivial partition $[1,m]= P_1 \sqcup \ldots \sqcup P_j$
such 
\begin{enumerate}
\item $P_j \cap [1,n] \neq \varnothing$ for all $j$ and 
\item $\wt{B}$ becomes a block-diagonal matrix with respect to this partition (with blocks of rectangular size).
\end{enumerate}
It is well known that, if the exchange matrix of one seed of $\AA$ is indecomposable, then this is true for any other seed of $\AA$. 
Analogously to \cite[Lemma 2.3]{ASS}, one proves:
\ble{indec} If the exchange matrix on one {\em{(}}and thus any{\em{)}} seed of $\AA$ is indecomposable, then an automorphism $\phi \in \Aut(\AA)$ 
is a modular Poisson automorphism if and only if there is a labeled seed $\Sig=(y_1, \ldots, y_m, \wt{B})$ of $\AA$ satisfying the following condition:

Either $((\phi(y_1), \ldots, \phi(y_m)), \wt{B})$ or $((\phi(y_1), \ldots, \phi(y_m)), -\wt{B})$ is a labeled seed of $\AA$ and the corresponding compatible 
$\La$-matrix is $\La(\Sig)$.
\ele
Once again, it is easy to check that if one seed of $\AA$ satisfies the condition in \leref{indec}, then every labeled seed of $\AA$ will satisfy it too. 
Denote by $\MAut (\AA, \{.,.\})$ the group of modular Poisson automorphisms of $(\AA, \{.,.\})$. The subgroup of those $\phi \in \Aut(\AA)$ satisfying 
the condition that $((\phi(y_1), \ldots, \phi(y_m)), \wt{B})$ is a labeled seed of $\AA$, will be called the group of {\em{strong modular Poisson automorphisms}}
of $(\AA, \{.,.\})$ in analogy with the terminology of \cite{FZ2} for cluster algebra automorphisms without the presence of Poisson structure.  

Finally, note that $\MAut (\AA, \{.,.\})$ normalizes $T_\AA$. Thus, we have the subgroup
\[
T_\AA \ltimes \MAut (\AA, \{.,.\}) \subseteq \Aut (\AA, \{.,.\}).
\]
The group $\MAut (\AA, \{.,.\})$ is of combinatorial nature and one can hope that it can be fully determined from the exchange pattern 
of the cluster algebra $\AA$. It has been described for acyclic cluster algebras \cite{ASS} and for cluster algebras of finite mutation type \cite{G1,BS,BY}
in terms of the transjective component of the Auslander--Reiten quiver of the corresponding cluster category and the mapping class group of the 
associated surface, respectively. 
The group $T_\AA$ is a torus that is described explicitly. However, the full automorphism group $\Aut (\AA, \{.,.\})$ is extremely hard to determine explicitly.
This rises to the following natural question: 
\bprob{prob} Let $\AA$ be a cluster algebra of geometric type over a field $\KK$, equipped with a Gekhtman--Shapiro--Vainshtein Poisson structure $\{.,.\}$.
Determine how much bigger the full Poisson automorphism group $\Aut (\AA, \{.,.\})$ is compared to its canonical subgroup
\[
T_\AA \ltimes \MAut (\AA, \{.,.\}).
\]
\eprob
In the next section we present a solution of this problem for the important case of the cluster algebra structures on all open Schubert cells 
in full flag varieties of complex simple Lie groups. We expect that \thref{4} will be helpful to fully resolve the above problem for all
$\Nset$-graded connected cluster algebras.
\sectionnew{Applications to Poisson automorphism of Schubert cells}
\lb{App-Sc}
\subsection{Standard Poisson structures on flag varieties} 
\label{4.1}
Let $\KK$ be a field of characteristic 0 and $G$ be a split, connected, simply connected, semisimple algebraic group over $\KK$. Let
$\g$ be its Lie algebra and $W$ its Weyl group. Denote by $\{\al_1, \ldots, \al_r\}$ 
the set of its simple roots and by $\{s_1, \ldots, s_r\} \subset W$ the corresponding set of simple reflections. 

Let $B_\pm$ be a pair of opposite Borel subgroups of $G$, 
$H: = B_+ \cap B_-$ be the corresponding maximal torus, and $U_\pm$ be the unipotent radicals of $B_\pm$. Set $\b_\pm : = \Lie (B_\pm)$, 
$\n_\pm := \Lie(U_\pm)$, and $\h := \Lie H$. 
Let $\{e_i, f_i\}$ be a set of Chevalley generators of $\g$ such that $\{e_i\}$ and $\{f_i\}$ generate $\n_+$ and $\n_-$, respectively. 
Set $h_i := [e_i, f_i]$.
Consider the representatives of $s_i$ in the normalizer $N_G(H)$ of $H$ in $G$
\[
\dot{s}_i := \exp(f_i) \exp (-e_i) \exp (f_i). 
\]
They are extended (in a unique way) to Tits' representatives of the Weyl group elements $w \in W$ in $N_G(H)$ by setting $\dot{w} := \dot {v} \dot{s}_i$ if 
$w = v s_i$ and $\ell(w) = \ell(v) + 1$ where $\ell \colon W \to \Nset$ is the length function. Denote by $\De_+$ the set of positive roots of $\g$. 
For $w \in W$ and $i \in [1,r]$ such that $w(\al_i) \in \De_+$, denote the root vectors 
\begin{equation}
\label{roots}
e_{w(\al_i)} := \Ad_{\dot{w}}(e_{\al_i}) \quad \mbox{and} \quad f_{w(\al_i)} := \Ad_{\dot{w}}(f_{\al_i}).
\end{equation}
Let $\lcor ., . \rcor$ be the invariant bilinear form on $\g$, normalized by $\|\al_i\|^2 =2$ for short roots $\al_i$.

The standard $r$-matrix for $\g$ is the element
\begin{equation}
\label{stand-r}
r_{\st}:= \sum_{\be \in \De_+} \frac{\| \be \|^2}{2} e_\be \wedge f_\be \in \wedge^2 \g.
\end{equation}
It gives rise to the Poisson bivector field
\begin{equation}
\label{pi}
\pi_{\st} := - \chi(r) \in \Ga(G/B_+, \wedge^2 T (G/B_+)),
\end{equation}
called the {\em{standard Poisson structure of the full flag variety}} $G/B_+$. Here, $\chi \colon \g \to \Ga(G/B_+, T (G/B_+))$ denotes 
the infinitesimal action of $\g$ associated to the left action of $G$ on $G/B_+$, and its extension to $\wedge^\bullet \g \to \Ga(G/B_+, \wedge^\bullet T (G/B_+))$.
Alternatively, the Poisson structure $\pi$ can be defined as the push-forward under the projection $G \mt G/B_+$ of the 
standard (Sklyanin) Poisson structure on the group $G$. The Schubert cell partition of $G/B_+$
\[
G/B_+ = \bigsqcup_{w \in W} B_+ w B_+/B_+
\]
is a decomposition into a Poisson submanifolds. 
Denote by $w_\circ$ the longest element of the Weyl group $W$ and consider the open Schubert cell of $G/B_+$
\[
B_+ w_\circ B_+/B_+ \cong U_+
\]
where the isomorphism (of affine spaces) is given by $u_+ \in U_+ \mt u_+ w_\circ B_+/B_+$. Denote by $\pi_{\st}$ the induced Poisson structure on 
$U_+$ and by $\{.,.\}_{\st}$ the corresponding Poisson bracket on $\KK[U_+]$. (Note that $\pi_{\st}$ is not the restriction of the standard 
Poisson structure of $G$ to $U_+$; even more, $U_+$ is not a Poisson submanifold of $G$ with respect to it.) 
Elek and Lu proved \cite[Theorem 6.1]{EL} that the Poisson structure $\pi_{\st}$ is defined over $\Zset$, see the next subsection for details.

Let $U_q(\g)$ be the corresponding quantized universal enveloping 
algebra \cite{L} and $U_q(\n_-)$ be its negative part. The Poisson bracket $\{.,.\}_{\st}$ is precisely the induced Poisson bracket for the 
specialization of an integral form of $U_q(\n_-)_{\opp}$ at $q=1$ (identified with the nonrestricted integral form on the 
quantized coordinate ring of the open Schubert cell $B_+ w_\circ B_+/B_+$), \cite[Lemma 4.3]{Y-strat}. 
Here $(.)_{\opp}$ refers to the opposite algebra structure.

The automorphism group $\Aut(\Ga)$ of the Dynkin graph $\Ga$ of $\g$ acts on $G$, $G/B_+$ and $U_+$ 
(via its action on the fixed choice of Chevalley generators of $\g$). It is 
clear that this action is Poisson with respect to $\pi_{\st}$. The definition of $\pi_{\st}$ also implies that 
the conjugation action of $H$ on $(U_+, \{.,.\}_{\st})$ is Poisson. 
Thus, we have the embedding
\begin{equation}
\label{em}
(H/Z_G)  \ltimes \Aut(\Ga)  \hra \Aut(\KK[U_+], \{.,.\}_{\st})
\end{equation}
where the semidirect product is formed with respect to the action of $\Aut(\Ga)$ on $H$ permuting the elements $\{h_i\}$
and $Z_G$ denotes the center of $G$.
By abuse of notation, we will identify $(H/Z_G) \ltimes \Aut(\Ga)$ with its image in $\Aut(\KK[U_+], \{.,.\}_{\st})$.

\bth{5} For all fields $\KK$ of characteristic $0$ and split, connected, simply connected, semisimple algebraic groups $G$
which do not have $SL_2$ factors, the group of Poisson automorphisms of the coordinate ring 
of the open Schubert cell $(B_+ w_\circ B_+/B_+, \pi_{\st}) \cong (U_+, \pi_{\st})$ is given by 
\[
\Aut(\KK[U_+], \{.,.\}_{\st}) \cong (H/Z_G) \ltimes \Aut(\Ga).
\]
\eth
For all semisimple Lie groups $G$, Berenstein, Fomin, and Zelevinsky constructed in \cite{BFZ} an upper cluster algebra structure on $\KK[U_+]$ 
with a Gekhtman--Shapiro--Vainshtein Poisson structure equal to $\{.,.\}_{\st}$. In the case when $G$ is simply laced, Geiss, Leclerc and Schr\"oer proved 
in \cite{GLS0} that $\KK[U_+]$ coincides with the corresponding cluster algebra $\AA(U_+)$. For all semisimple Lie algebras this was done
by Goodearl and the second author in \cite{GY-P}. \thref{5} solves \probref{prob} for this class of cluster algebras. 
It also has the following immediate corollary for the Poisson modular automorphism groups of the cluster algebras $\AA(U_+)$:
\bco{5b} For all fields $\KK$ of characteristic $0$ and split, connected, simply connected, semisimple algebraic groups $G$
which do not have $SL_2$ factors, the Poisson modular 
automorphism group of the cluster algebra structure $(\AA(U_+), \{.,.\}_{\st})$ on $\KK[U_+]$ is isomorphic to
$\Aut(\Ga)$.
\eco 
\thref{5} and \coref{5b} are not valid in the case when the semisimple algebraic group $G$ has $SL_2$ factors because the Poisson 
structure $\pi_{\st}$ vanishes in the case of $G=SL_2$. This leads to a slightly larger automorphism group in the case of one factor and to pathological 
problems in the case of multiple factors. In the latter case the Poisson algebra in question is a tensor product of two algebras, one of which is a polynomial 
algebra with a trivial Poisson bracket.

Since the category of Poisson algebras in \thref{5} is closed under tensor products, we have the following application
of \thref{5} that solves the isomorphism problem for this family of Poisson algebras. For a split, connected, simply connected, semisimple 
algebraic group $G$, we will denote by $U_+(G)$ the unipotent radical of a Borel subgroup of $G$ and by 
$ \{.,.\}_{\st}$ the corresponding standard Poisson structure on $\KK[U_+]$. The Poisson algebra $(\KK[U_+(G)], \{.,.\}_{\st})$
is independent (up to isomorphism) on the choice of Borel subgroup of $G$.

\bth{5c} Let $G_1$ and $G_2$ be split, connected, simply connected, semisimple algebraic groups
over a field $\KK$ of characteristic $0$. The Poisson algebras 
\[
(\KK[U_+(G_1)], \{.,.\}_{\st}) \quad \mbox{and} \quad (\KK[U_+(G_2)], \{.,.\}_{\st})
\]
are isomorphic, if and only if $G_1$ and $G_2$ are isomorphic.
\eth
\thref{5c} does not assume that the semisimple groups do not have $SL_2$ factors, like \thref{5}. \thref{5} is proved in \S \ref{4.2}-\ref{4.5} 
and \thref{5c} in \S \ref{4.6}.
\subsection{Generators for the Poisson algebra $(\KK[U_+], \{.,.\}_{\st})$}
\label{4.2}
Fix a reduced expression of the longest Weyl group element
\begin{equation}
\label{reduced}
w_\circ = s_{i_1} \ldots s_{i_N}
\end{equation}
where $N:= \dim \n_+$. For $k \in [1,N]$, set
\[
(w_\circ)_{< k} := s_{i_1} \ldots s_{i_{k-1}} \in W.
\]
The positive roots of $\g$ are expressed as
\[
\be_k =(w_\circ)_{< k}(\al_{i_k}), \quad k \in [1,N]
\]
and, using \eqref{roots}, we have the parametrization 
\[
U_+ = \{ \exp(x_{\be_1} e_{\be_1}) \ldots \exp(x_{\be_N} e_{\be_N}) \mid x_{\be_k} \in \KK \}.
\]
It gives rise to the isomorphism
\begin{equation}
\label{coord-U}
\KK[U_+] \cong \KK[x_\be, \be \in \De_+].
\end{equation}
These are the coordinates in which the Poisson structure $\{.,.\}_{\st}$ has integer coefficients, \cite[Theorem 6.1]{EL}.
It also satisfies the semiclassical analog of the Levendorskii--Soibelman straightening law
\begin{equation}
\label{LS}
\{x_{\be_j}, x_{\be_k}\}_{\st} = \lcor \be_j, \be_k \rcor x_{\be_j} x_{\be_k} + p_{jk} \quad
\mbox{where} \; \; p_{jk} \in \KK[x_{\be_{j+1}}, \ldots, x_{\be_{k-1}} ]
\end{equation}
for all $j<k \in [1,N]$, \cite[Proposition 5.12]{EL}. (We use the opposite Poisson structure to that in \cite{EL}, resulting
in a sign difference.)

Denote by $Q_+ = \Nset \al_1 + \ldots + \Nset \al_r$ the positive part of the root lattice of $\g$ and by $\htt \colon Q_+ \to \Nset$ the 
principal grading 
\[
\htt(n_1 \al_1 + \cdots + n_r \al_r) := n_1 + \cdots + n_r.
\]
The algebra $\KK[U_+]$ is $Q_+$-graded by setting $\deg x_\be := \be$ for $\be \in \De_+$. 
(One can equivalently define the grading by using the characters of the conjugation action of $H$ on $\KK[U_+]$.)
By \eqref{pi}, the Poisson structure $\{.,.\}_{\st}$ is graded of degree 0.

\bre{change-red} Since $\KK[U_+]_{\al_i} = \KK x_{\al_i}$ for $i \in [1,r]$, changing the reduced expression of $w_\circ$ only 
results in rescaling of the coordinate functions $x_{\al_i}$.
\ere

\bpr{gen} The Poisson algebra $(\KK[U_+], \{.,.\}_{\st})$ is generated by $x_{\al_i}$ for $i \in [1,r]$.
\epr
\begin{proof} The ordering  $\be_1, \ldots, \be_N$  of the positive roots of $\g$ is convex in the sense that, if $\al, \be, \al+\be \in \De_+$, then 
$\al + \be$ appears between $\al$ and $\be$. 

Denote by $\PP$ the Poisson subalgebra of $\KK[U_+]$ generated by $x_{\al_i}$ for $i \in [1,r]$.
We prove by induction on $n \in \Zset_{\geq 1}$ that $x_\be \in \PP$ for $\be \in \De_+$ with $\htt(\be)=n$. This is clear for $n=1$. Assume its validity for some 
$n \in \Zset_{\geq 1}$ and consider $\ga \in \De_+$ with $\htt(\ga) = n+1$. Then $\ga = \be + \al_i$ for some $\be \in \De_+$, $\htt(\be)=n$. 
Since  $\{.,.\}_{\st}$ is graded of degree 0 with respect to the $Q_+$-grading of $\KK[U_+]$, 
\[
\{x_{\al_i}, x_\be \}_{\st} = a x_\ga + p 
\]
where $a \in \KK$ and $p$ is a polynomial in $x_{\be'}$ with $\be' \in \De_+$, $\htt(\be') \leq n$. The inductive assumption implies that
$p, \{x_{\al_i}, x_\be \}_{\st} \in \PP$, so it suffices 
to prove that 
\begin{equation}
\label{aneq}
a \neq 0. 
\end{equation}
Let $j,k \in [1,N]$ be such that $\be_j = \al_i$ and $\be_k = \be$. 
By the convexity of the ordering of the roots $\be_1, \ldots, \be_N$, $\ga = \be_l$ for some $l$ between $j$ and $k$.
Denote $m = \min \{j,k\}$ and $M := \max\{j, k\}$.  
A formula for the coefficient $a$ was obtained in \cite[Theorem 4.10, Eq. (41)]{EL}. Eq. (41) in \cite{EL}
expresses $a$ as a product of 1's and a nonzero rational number coming from the $\al_{i_l}$-string through the root 
\[
s_{i_{l-1}} \ldots s_{i_{m-1}}(\al_{i_m}) - \al_{i_l} =  s_{i_{l-1}} \ldots s_{i_{M-1}}(\al_{i_M}) \in \De_+.
\]
The last identity follows from the equality $\ga = \be + \al_i$ and the definitions of $l,m$ and $M$.
This completes the induction step.
\end{proof}
\subsection{Poisson clusters of the algebra $(\KK[U_+], \{.,.\}_{\st})$}
\label{4.3}
Denote by $\{\vpi_1, \ldots, \vpi_r\}$ the fundamental weights  of $G$.
Given a dominant integral weight $\la$ of $G$ and $u, v \in W$, one defines the generalized minors
\[
\De_{u \la, v \la} \in \KK[G] 
\]
as follows. Let $L(\la)$ denote the irreducible highest weight $G$-module with highest weight $\la$. Fix  a highest weight vector $b_\la$
of $L(\la)$ and a vector $\xi_\la$ in the dual weight space, normalized by $\lcor \xi_\la, b_\la \rcor =1$. Define
\[
\De_{u \la, v \la} (g) := \lcor \xi_\la, \dot{u}^{-1} g \dot{v} b_\la \rcor, \quad g \in G.
\]
As in the previous subsection, we fix a reduced expression \eqref{reduced} of $w_\circ$ and define 
\begin{equation}
\label{normal-De}
\De_k:= \De_{(w_\circ)_{<k} \vpi_k, w_\circ \vpi_k} \in \KK[U_+] \quad \mbox{for} \; \; k \in [1,N].
\end{equation}
(The minors are considered as elements of $\KK[U_+]$ by restriction.) Clearly, for all $i \in [1,r]$, 
\[
\De_{\vpi_i, w_\circ \vpi_i} = \De_{k} \quad \mbox{where} \; \; k := \min \{l \in [1,N] \mid i_l = i \}.
\]
\bpr{Poisson-cl} For each reduced expression of $w_\circ$, the elements $\De_1, \ldots, \De_N \in \KK[U_+]$ form a 
Poisson cluster of $(\KK[U_+], \{.,.\}_{\st})$ and satisfy 
\begin{equation}
\label{P-bracket}
\{ \De_j, \De_k \}_{\st} = - \lcor ((w_\circ)_{<j} + w_\circ) \vpi_{i_j}, ((w_\circ)_{< k} - w_\circ) \vpi_{i_k} \rcor \De_j \De_k, \quad \forall j<k \in [1,N].
\end{equation}
\epr
\begin{proof} Denote the successor function 
for the reduced expression \eqref{reduced}
\[
s \colon [1,N] \to [1,N] \sqcup \{ \infty \}
\]
given by
\[
s(k) = 
\begin{cases}
\min \{ l> k \mid i_l = i_k \}, & \mbox{if $\exists$ $l>k$ such that $i_l =i_k$} 
\\
\infty, &\mbox{otherwise}.
\end{cases}
\]

The fact that $\De_1, \ldots, \De_N$ form a polynomial subring of $\KK[U_+]$ and the inclusion 
\[
\KK[U_+] \subset \KK[\De_1^{\pm 1}, \ldots, \De_N^{\pm1}]
\]
follow from the expansion
\begin{equation}
\label{expan}
\De_k = \De_{s(k)} x_{\be_k} + p_k \quad \mbox{for some} \; \; p_k \in \KK[x_{\be_{k+1}}, \ldots, x_{\be_N}], \quad \forall k \in [1,N]
\end{equation}
where $\De_{\infty}:= 1$ and $p_N =0$. It is proved analogously to its quantum counterpart \cite[Proposition 3.3]{GeY}. 

The Poisson structure $\{.,.\}_{\st}$ on the open Schubert cell $B_+ w_\circ B_+/B_+$ equals the canonical Poisson
structure on the specialization at $q=1$ of the nonrestricted integral form of the quantized coordinate ring of $B_+ w_\circ B_+/B_+$, 
see \cite[pp. 274-276]{Y-strat}. The same arguments show that the minors $\De_1, \ldots, \De_N$ are equal to
the specializations of the corresponding set of quantum minors in \cite[Theorem 9.5]{GY-qcl} for the case of $w=w_\circ$. 
The formula for Poisson brackets \eqref{P-bracket} is an immediate consequence of the commutation relations between 
the quantum minors in \cite[Eq. (9.29)]{GY-qcl}, with a sign contribution from the antihomomorphisms
in \cite[Theorem 9.2]{GY-qcl}.
\end{proof}
For each $i \in [1,r]$, $- w_\ci(\al_i)$ is a positive simple root of $\g$. Denote the involution $\tau$ of $[1,r]$, given by
\begin{equation}
\label{tau}
\al_{\tau(i)}:= - w_\ci(\al_i).
\end{equation}
Let $O_1 \subset [1,r]$ be the set of its fixed points and $O_2$ be a set or representatives of its 2-element orbits.

\bpr{c-Poisson-cl} For each reduced expression of $w_\circ$, the Poisson center of the Laurent polynomial ring 
of the Poisson cluster in \prref{Poisson-cl}
is given by
\[
\ZZ(\KK[\De_1^{\pm 1}, \ldots, \De_N^{\pm 1}]) = \KK [ \De_{\vpi_i, w_\circ \vpi_i}^{\pm 1}, (\De_{\vpi_l, w_\circ \vpi_l} \De_{\vpi_{\tau(l)}, w_\circ \vpi_{\tau(l)}})^{\pm 1},
i \in O_1, l \in O_2].
\]
\epr
\begin{proof} It follows from \eqref{P-bracket} that the Poisson center of the Laurent polynomial ring 
$\ZZ(\KK[\De_1^{\pm 1}, \ldots, \De_N^{\pm 1}])$ is given in terms of the radical of the form $\Om$ on $\Zset^N$ 
\[
\Om(\de_j, \de_k)  := - \lcor ((w_\circ)_{<j} + w_\circ) \vpi_{i_j}, ((w_\circ)_{< k} - w_\circ) \vpi_{i_k} \rcor, \quad j < k \in [1,N]
\]
by \eqref{cent}. (Recall the notation \eqref{stand-basis}.)
By \cite[Eq. (9.29)]{GY-qcl} the radical of the same form computes the center of the quantum torus generated by the 
quantum counterparts of this sequence of minors. The center of that quantum torus was described in \cite[Lemma 4.7]{Y-AD}.
The proposition follows from it and \eqref{normal-De}.
\end{proof}
\subsection{Proof of \thref{5} in the unipotent case} 
\label{4.4}
Consider the principal grading of the algebra $\KK[U_+]$, defined by 
\begin{equation}
\label{prin}
\deg x_\be = \htt(\be) \quad \mbox{for} \; \; \be \in \De_+.
\end{equation}
The grading is connected, i.e., $\KK[U_+]_0 = \KK$. 

We prove \thref{5} in two steps. In this subsection we prove it in the case of unipotent automorphisms, using the 
rigidity theorem from Sect. \ref{Rig} and the embedding result for unipotent automorphisms on the basis of Poisson clusters in \prref{aut-emb}. 
In the next subsection we reduce the general case to the unipotent one (with respect to the 
principal grading of $\KK[U_+]$.

For the needs of the proofs in the next subsection, we derive a stronger result for unipotent automorphisms
concerning all specializations of the $Q_+$-grading on $\KK[U_+]$ to $\Nset$-gradings. Denote the set of 
dominant integral coweights of $G$:
\[
P_+\spcheck = \{ n_1 \vpi_i\spcheck + \cdots + n_r \vpi_r\spcheck \mid n_1, \ldots, n_r \in \Nset \}
\]
where $\vpi_i\spcheck, \ldots, \vpi_r\spcheck$ are the fundamental coweights of $G$.
Each $\la \in P_+\spcheck$ gives rise to the specialization of the $Q_+$-grading of $\KK[U_+]$, defined by
\begin{equation}
\label{la-grad}
\deg x_{\be}:= \lcor \la, \be \rcor.
\end{equation}

\bpr{unip-wci} Let $\KK$ be a field of characteristic $0$ and $G$ be a split, connected, simply connected, semisimple algebraic group.
For all dominant integral coweights $\la$ of $G$, the group of unipotent Poisson automorphisms of the coordinate ring of 
the open Schubert cell $(B_+ w_\circ B_+/B_+, \pi_{\st}) \cong (U_+, \pi_{\st})$ with respect to the grading \eqref{la-grad} is trivial
\[
\UAut(\KK[U_+], \{.,.\}_{\st}) = \{\id\}.
\]
\epr
For this proposition it is not necessary to require that $G$ has no $SL_2$ factors.
\begin{proof} Since the scheme theoretic intersection of opposite Schubert varieties is reduced \cite[Theorem 3.5]{R} and Schubert varieties are linearly defined 
\cite[Theorem 3(i)]{KR}, the vanishing ideal of the irreducible subvariety
\[
U_+ \cap \ol{U_- s_i B_+} w_\circ \cong B_+ w_\circ B_+/B_+ \cap \ol{B_- s_i B_+/B_+} 
\]
of $U_+$
is the principal ideal $(\De_{\vpi_i, w_\circ \vpi_i})$. Because $\KK[U_+]$ is a polynomial ring, the elements $\De_{\vpi_i, w_\circ \vpi_i} \in \KK[U_+]$ 
are prime. Clearly, they are not associates of $x_{\al_l}$ for any $l \in [1,r]$.

Let $\phi \in \UAut(\KK[U_+], \{.,.\}_{\st})$ with respect to the grading \eqref{la-grad}. Fix $i \in [1,r]$. First we show that 
\begin{multline}
\label{phi-x-al}
\phi(x_{\al_i}) = (1 + f_i) x_{\al_i} \quad \mbox{for some} \\
f_i \in \KK [ \De_{\vpi_j, w_\circ \vpi_j}, \De_{\vpi_l, w_\circ \vpi_l} \De_{\vpi_{\tau(l)}, w_\circ \vpi_{\tau(l)}},
j \in O_1, l \in O_2]_{\geq 1}.
\end{multline}
Recall the setting of \prref{c-Poisson-cl} and the definition \eqref{tau} of the involution $\tau$ of the Dynkin graph $\Ga$. 
Keeping in mind \reref{change-red}, we choose a reduced expression \eqref{reduced} of $w_\circ$
such that $i_N =\tau(i)$. So, $\be_N = w_\ci( s_{\tau(i)} \al_{\tau(i)}) = w_\ci(- \al_{\tau(i)}) = \al_i$. 
By \prref{Poisson-cl}, the $N$-tuple $(\De_1, \ldots, \De_N)$ is a Poisson cluster of $(\KK[U_+], \{.,.\}_{\st})$. 
In addition, it follows from \eqref{expan} that 
\[
\De_N = x_{\be_N} = x_{\al_i}.
\]

We apply \prref{aut-emb} to obtain a unipotent bi-integral automorphism $\iota(\phi)$ of $(\KK[\De_1^{\pm 1}, \ldots, \De_N^{\pm 1}], \{.,.\}_{\st})$.
The rigidity result in \thref{1} and the classification of the Poisson center of this cluster from \prref{c-Poisson-cl} give that 
\begin{multline*}
\phi(x_{\al_i}) = \phi(\De_N) = \iota(\phi)(\De_N) = (1 + f_i) x_{\al_i} \quad \mbox{with} \\
f_i \in \KK [ (\De_{\vpi_j, w_\circ \vpi_j})^{\pm 1}, (\De_{\vpi_l, w_\circ \vpi_l} \De_{\vpi_{\tau(l)}, w_\circ \vpi_{\tau(l)}})^{\pm 1},
j \in O_1, l \in O_2]_{\geq 1}.
\end{multline*}
However, $\De_{\vpi_j, w_\circ \vpi_j}$ is a prime element of $\KK[U_+]$ for all $j \in [1,r]$ which is not a multiple of $x_{\al_i}$. Hence,
\[
f_i \in \KK [ \De_{\vpi_j, w_\circ \vpi_j}, \De_{\vpi_l, w_\circ \vpi_l} \De_{\vpi_{\tau(l)}, w_\circ \vpi_{\tau(l)}},
j \in O_1, l \in O_2]_{\geq 1}
\]
which proves \eqref{phi-x-al}. Since $f_i$ is in the Poisson center of $(\KK[U_+], \{.,.\}_{\st})$ and $\KK[U_+]$ is 
generated as a Poisson algebra by $x_{\al_i}$, $i \in [1,r]$ (\prref{gen}), 
\[
\phi(\De_{\vpi_i, w_\circ \vpi_i}) = \De_{\vpi_i, w_\circ \vpi_i} \prod_{l=1}^r (1+ f_l)^{m_{il}}, \quad \forall i \in [1,r]
\]
where $m_{il} \in \Nset$ are the integers given by
\[
\vpi_i- w_\circ \vpi_i = \sum_l m_{il} \al_l.
\]
Because $\De_{\vpi_i, w_\circ \vpi_i}$ is a prime element of $\KK[U_+]$ and 
$\phi$ is an automorphism of $\KK[U_+]$, 
\[
\prod_{l=1}^r (1+ f_l)^{m_{il}}=1, \quad \forall i \in [1,r]. 
\]
We also have $f_l \in \KK[U_+]_{\geq 1}$. 
Taking into account the fact that for each $l \in [1,r]$
there exists $i \in [1,r]$ such that $m_{il} \neq 0$, we obtain that all $f_l=0$. This completes the proof of the theorem.
\end{proof}

\subsection{Reduction of \thref{5} to the case of unipotent automorphisms}
\label{4.5}
\bpr{reduct} Let $\KK$ be a field of characteristic $0$ and $G$ be a split, connected, simply connected, semisimple algebraic group
which does not have $SL_2$ factors. Then, for every $\phi \in \Aut(\KK[U_+], \{.,.\}_{\st})$ there exists  $\psi \in (H/Z_G) \ltimes \Aut(\Ga)$
such that $\phi \psi^{-1} \in \UAut(\KK[U_+], \{.,.\}_{\st})$ with respect to the principal grading.
\epr

Unless otherwise noted, all results in this section refer to the principal grading \eqref{prin} of $\KK[U_+]$.
First we show that every automorphism of $(\KK[U_+], \{.,.\}_{\st})$ is increasing with respect to this grading:  
\ble{1} In the setting of \thref{5}, every automorphism $\phi$ of $(\KK[U_+], \{.,.\}_{\st})$ satisfies
\[
\phi(f) \in \KK[U_+]_{\geq n}, \quad \mbox{for all} \; \; f \in \KK[U_+]_n, n \in \Nset.
\]
\ele
\begin{proof} Since $\KK[U_+]$ is generated as a Poisson algebra by $x_{\al_i}$, $i \in [1,r]$ (\prref{gen}) and $\{.,.\}_{\st}$ has degree 0, 
it is sufficient to prove that 
\begin{equation}
\label{phi-x}
\phi(x_{\al_i}) \in \KK[U_+]_{\geq 1}.
\end{equation}
The definition of $x_{\al_i}$ does not depend, up to a scalar, on the choice of reduced expression of $w_\circ$ (\reref{change-red}).
Because $G$ does not have $SL_2$ factors, there exists a simple root $\al_l$ such that $\lcor \al_i, \al_l \rcor \neq 0$. Recall the definition \eqref{tau} of the 
involution $\tau$ of the Dynkin graph $\Ga$. Choose a reduced expression of $w_\circ$ ending with $s_{\tau(l)} s_{\tau(i)} \in W$, i.e., an expression \eqref{reduced} 
such that $i_{N-1} = \tau(l)$ and $i_N =\tau(i)$. Thus, $\be_N = w_\ci( s_{\tau(i)} \al_{\tau(i)}) = w_\ci(- \al_{\tau(i)}) = \al_i$, 
$\be_{N-1} = w_\ci( s_{\tau(i)} s_{\tau(l)} (\al_{\tau(l)}))= w_\ci(- s_{\tau(i)} \al_{\tau(l)}) = s_i \al_l$ and
\[
\lcor \be_{N-1}, \be_N \rcor = \lcor s_i (\al_l), \al_i \rcor = - \lcor \al_l, \al_i \rcor \neq 0.
\]
By \eqref{LS},
\[
\{x_{\be_{N-1}}, x_{\al_i} \}_{\st} = - \lcor \al_l, \al_i \rcor x_{\be_{N-1}} x_{\al_i}. 
\] 
Let $\phi(x_{\be_{N-1}}) = f + g$ with $f \in \KK[U_+]_n$, $f \neq 0$ and $g \in \KK[U_+]_{>n}$ for some $n \in \Nset$. Assume that \eqref{phi-x} does not hold.
Then $\phi(x_{\al_i}) - a \in \KK[U_+]_{\geq 1}$ for some $a \in \KK^*$. The term of minimal degree in the RHS of the equality
\[
\{ \phi(x_{\be_{N-1}}), \phi(x_{\al_i}) \}_{\st} = - \lcor \al_l, \al_i \rcor \phi(x_{\be_{N-1}}) \phi(x_{\al_i})
\]
has degree $n$ and equals $- \lcor \al_l, \al_i \rcor af \neq 0$, while the term of degree $n$ in the LHS equals $\{f,a \}_{\st} =0$. 
This is a contradiction, which proves \eqref{phi-x} and the lemma. The last step of the proof is a Poisson analog of the argument 
of \cite[Proposition 3.2]{LaLe}. 
\end{proof}

We will call an automorphism $\phi$ of $(\KK[U_+], \{.,.\}_{\st})$ {\em{linear}} if  
\[
\phi(\KK[U_+]_n) = \KK[U_+]_n, \; \; \forall n \in \Nset, \quad \mbox{i.e.}, \; \; \phi(x_{\al_i}) \in \KK[U_+]_1, \; \; \forall i \in [1,r]
\]
with respect to the principal grading.
The following result classifies those automorphisms.
\bpr{2} In the setting of \thref{5}, the group of linear automorphisms of $(\KK[U_+], \{.,.\}_{\st})$ is isomorphic to 
$(H/Z_G) \ltimes \Aut(\Ga)$.
\epr
For the proof of this fact we consider the presentation of $(\KK[U_+], \{.,.\}_{\st})$ as a factor of a free 
Poisson algebra. Recall that $\KK[U_+]$ is generated by $x_{\al_i}$ as a Poisson algebra (\prref{gen}). 
Denote by $(\FF_G, \{.,.\})$ the free Poisson algebra on $x_{\al_i}$ and by $\II_G$ the kernel of the 
canonical projection of Poisson algebras $(\FF_G, \{.,.\}) \to (\KK[U_+], \{.,.\}_{\st})$. Consider the $Q_+$-grading of $\FF_G$, given by $\deg x_{\al_i} := \al_i$,
$i \in [1,r]$.
The projection is $Q_+$-graded. Thus, $\II_G$ is homogeneous with respect to the $Q_+$-grading and its principal specialization.

Denote by $(c_{ij})$ the Cartan matrix of $G$.  
\ble{3} The degree 3 component of the ideal $\II_G$ is spanned by 
\begin{equation}
\label{I-1}
x_{\al_l} \{ x_{\al_i}, x_{\al_j} \}, \{x_{\al_l}, \{ x_{\al_i}, x_{\al_j} \} \}
\end{equation}
for $l,i, j \in [1,r]$ such that $c_{ij}=0$ and 
\begin{equation}
\label{I-2}
\{ x_{\al_i}, \{ x_{\al_i}, x_{\al_j} \} \} - 
\lcor \al_i, \al_j \rcor^2 x_{\al_i}^2 x_{\al_j} 
\end{equation}
for $i, j \in [1,r]$ such that $c_{ij} =-1$.
\ele
\begin{proof} Fix $i \neq j \in [1,r]$. Consider a reduced expression \eqref{reduced} of $w_\circ$ such that 
$i_1 =i$ and $i_2 =j$.

If $c_{ij} =0$, then $\be_2 = \al_j$ and by \eqref{LS},
\[
\{ x_{\al_i}, x_{\al_j} \}_{\st} = \lcor \al_i, \al_j \rcor x_{\al_i} x_{\al_j} = 0
\]
in $\KK[U_+]$. This implies that the elements in \eqref{I-1} belong to $\II_G$. 

If $c_{ij} =-1$, then $\be_2 = s_i (\al_j) = \al_i + \al_j$ and by \eqref{LS},
\begin{equation}
\label{1-P}
\{x_{\al_i}, x_{\al_i + \al_j} \}_{\st} = \lcor \al_i, s_i (\al_j) \rcor x_{\al_i} x_{\al_i + \al_j} = - \lcor \al_i, \al_j \rcor x_{\al_i} x_{\al_i + \al_j}
\end{equation}
in $\KK[U_+]$. Eq. \eqref{aneq} and the fact that $\KK[U_+]_{\al_i+\al_j}= \KK x_{\al_i} x_{\al_j} + \KK x_{\al_i + \al_j}$ imply
\[
\{x_{\al_i}, x_{\al_j} \}_{\st} = \lcor \al_i, \al_j \rcor x_{\al_i} x_{\al_j} + a x_{\al_i + \al_j}
\]
for some $a \in \KK^*$.
Expressing $x_{\al_i + \al_j}$ from here and substituting it in \eqref{1-P} shows that the elements of the type \eqref{I-2} are in the ideal $\II_G$.  
Finally, the property that $(\II_G)_\ga =0$ for the other $\ga \in \De_+$ with $\htt(\ga) =3$ follows from the isomorphism
$\KK[U_+] \cong \KK[x_{\be_1}, \ldots, x_{\be_N}]$ by comparing the dimensions of the degree 3 components of $\FF_G$ and $\KK[U_+]$. 
\end{proof}

The next lemma provides important restrictions on the possible form of linear automorphisms of $(\KK[U_+], \{.,.\}_{\st})$ 
coming from the degree 3 component of the ideal $\II_G$. 

\ble{4} Let $\phi$ be a linear isomorphism of $(\KK[U_+],\{.,.\}_{\st})$, given by 
\begin{equation}
\label{line}
\phi(x_{\al_i}) = \sum_{l=1}^r a_{il} x_{\al_l}
\end{equation}
for some $a_{il} \in \KK$. Denote the support of each row of this matrix
\[
\Supp(i) := \{ l \in [1,r] \mid a_{il} \neq 0 \}. 
\]
For all $i, j \in [1,r]$ such that $c_{ij} =-1$, we have the following:
\begin{enumerate}
\item[(i)] $\Supp(i) \cap \Supp(j) = \varnothing$, 
\item[(ii)] $|\Supp(i)|=1$, and 
\item[(iii)] for $k \in \Supp(i), m \in \Supp(j)$, $c_{km} =-1$.
\end{enumerate}
\ele

\begin{proof} Our strategy for the proof of this result is similar to that of the proof of \cite[Lemma 4.7]{Y-AD}.
Denote by $\wh{\phi}$ the automorphism of the free Poisson algebra $\FF_G$, given by the same formula \eqref{line} as $\phi$. Since 
$\phi$ is an automorphism of $\KK[U_+]$,
$\wh{\phi}(\II_G) = \II_G$. 

(i) Assume that $l \in \Supp(i) \cap \Supp(j)$. The degree $3\al_l$-component of the image of \eqref{I-2} under 
$\wh{\phi}$ is 
\[
- \lcor \al_i, \al_j \rcor^2 a_{il}^2 a_{jl} x_{\al_l}^3.
\]
It belongs to $\II_G$ because $\II_G$ is a graded ideal of $\FF_G$ with respect to its $Q_+$-grading.
This is in contradiction with \leref{3} which proves (i). 

(ii) Now assume that $k \neq l \in \Supp(i)$ and $m \in \Supp(j)$. By part (i), $k,l,m \in [1,r]$ are distinct. The degree 
$(\al_k + \al_l + \al_m)$-component of the image of \eqref{I-2} under $\wh{\phi}$ is a nonzero scalar multiple of 
\[
\{ x_{\al_k}, \{ x_{\al_l}, x_{\al_m} \} \} + \{ x_{\al_l}, \{ x_{\al_k}, x_{\al_m} \} \} - 2 \lcor \al_i, \al_j \rcor^2 x_{\al_k} x_{\al_l} x_{\al_m}.
\]
This expression belongs to $\II_G$ because $\II_G$ is $Q_+$-graded, but once again this contradicts \leref{3} since $(\II_G)_3$ 
does not contain elements of this form. This proves (ii).

(iii) The degree $(2 \al_k + \al_m)$-component of the image of \eqref{I-2} under $\wh{\phi}$ gives that
\[
\{x_{\al_k}, \{x_{\al_k}, x_{\al_m} \} \} - \lcor \al_k, \al_m \rcor^2 x_{\al_k}^2 x_{\al_m} \in \II_G.
\]
\leref{3} implies that $c_{km} = -1$. 
\end{proof}
\begin{proof}[Proof of \prref{2}] Eq. \eqref{em} is an embedding of $(H/Z_G) \ltimes \Aut(\Ga)$ in the group of linear automorphisms of 
$(\KK[U_+], \{.,.\}_{\st})$. We need to show that this embedding is an isomorphism. Two simple (positive) roots of $G$ will be
called adjacent if they are connected in the Dynkin graph $\Ga$.

Let $\phi$ be a linear automorphism 
of $(\KK[U_+], \{.,.\}_{\st})$. Applying \leref{4} to $\phi$ and $\phi^{-1}$, we obtain that there exists a permutation $\sigma \in S_r$
and $a_1, \ldots, a_r \in \KK^*$ such that:
\begin{enumerate}
\item[(i)] $\sigma$ maps the connected components of $\Ga$ of type $B_2$ and $G_2$ to connected components of type $B_2$ or $G_2$,
preserving the direction of arrows. Its restriction to the union of connected components of $\Ga$ that are not of type $B_2$ and $G_2$ is 
an isomorphism.

\item[(ii)] If $\al_i$ is the short simple root of a connected component of $\Ga$ of type $B_t$ for $t>2$, then 
\[
\phi(x_{\al_i}) = a_i x_{\al_{\sigma(i)}} + b_i x_{\sigma(j)}
\]   
where $\al_j$ is the second adjacent simple root to the only simple root of $\Ga$ adjacent to $\al_i$. (Note that by (i), the restriction
of $\sigma$ to a connected component of $\Ga$ of type $B_t$ for $t>2$ is an isomorphism.)

\item[(iii)] For all other simple roots $\al_i$ of $G$, 
\[
\phi(x_{\al_i}) = a_i x_{\al_{\sigma(i)}}.
\]
\end{enumerate}
Taking into account the property (iii), we see that $\sigma \in S_r$ cannot take a connected component of $\Ga$ of type $B_2$ 
to a connected component of type $G_2$ because in that case a restriction of $\phi$ will provide an isomorphism 
between the Poisson spaces $(\KK[U_+], \{.,.\}_{\st})$ for groups $G$ of type $B_2$ and $G_2$. This cannot happen because 
the corresponding unipotent radicals $U_+$ have different dimensions. Therefore, $\sigma$ is an automorphism of the Dynkin graph $\Ga$.
This implies that the map
\[
\phi_{\sigma, a} (x_{\al_i}) := a_i x_{\al_{\sigma(i)}}, \quad i \in [1,r]
\]
is a linear automorphism of $(\KK[U_+], \{.,.\}_{\st})$ coming from $(H/Z_G) \ltimes \Aut(\Ga)$. 
It remains to show that the scalars $b_i$ in property (ii) vanish, i.e., $\phi = \phi_{\sigma, a}$.  

Consider the dominant integral coweight 
\[
\la = n_1 \vpi_1\spcheck + \cdots + n_r \vpi_r\spcheck
\]
where $n_i=1$ if $\al_i$ is the short simple root of a connected component of $\Ga$ of type $B_t$ for $t>2$ and $n_i =2$
otherwise. It follows from properties (ii)--(iii) that $\phi \circ \phi_{\sigma, a}^{-1}$ is a unipotent automorphism of $(\KK[U_+], \{.,.\}_{\st})$
for the grading \eqref{la-grad}. \prref{unip-wci} implies that $\phi \circ \phi_{\sigma, a} =\id$ completing the proof of the proposition.
\end{proof}
\begin{proof}[Proof of \prref{reduct}.] Consider the principal grading of $\KK[U_+]$. 
Let $\phi$ be an automorphism of the Poisson algebra $\Aut (\KK[U_+], \{.,.\}_{\st})$. \leref{1} implies that for every $f \in \KK[U_+]_n$, $n \in \Nset$, 
there exists a unique  $\phi_0(f) \in \KK[U_+]_n$ such that 
\[
\phi(f) - \phi_0(f) \in \KK[U_+]_{> n}.
\]
This defines a map $\phi_0 \colon \KK[U_+] \to \KK[U_+]$ which is a linear automorphism of $(\KK[U_+], \{.,.\}_{\st})$ because 
$\{.,.\}_{\st}$ has degree 0 and $\phi_0 \circ (\phi^{-1})_0 = \id$.
It follows from \prref{2} that $\phi_0 \in (H/Z_G) \ltimes \Aut(\Ga)$, and clearly,
\[
\phi \phi_0^{-1} \in \UAut (\KK[U_+], \{.,.\}_{\st}).
\]
So, choosing $\psi :=\phi_0$ proves the desired property.
\end{proof}
\begin{proof}[Proof of \thref{5}.] \prref{reduct} implies that for every 
$\phi \in \Aut (\KK[U_+], \{.,.\}_{\st})$, there exists $\psi \in (H/Z_G) \ltimes \Aut(\Ga)$ such that $\phi \psi^{-1}$ 
is a unipotent automorphism of the Poisson algebra $(\KK[U_+], \{.,.\}_{\st})$ with respect to the principal grading of $\KK[U_+]$. 
\prref{2} implies that $\phi \psi^{-1} = \id$. Thus, $\phi = \psi$ which proves the theorem.
\end{proof}
\subsection{Proof of \thref{5c}}
\label{4.6}
Fix a Poisson algebra isomorphism 
\begin{equation}
\label{eta}
\eta \colon (\KK[U_+(G_1)], \{.,.\}_{\st}) \stackrel{\cong}{\lra}  (\KK[U_+(G_2)], \{.,.\}_{\st}).
\end{equation}
The map $\phi:= \eta \otimes \eta^{-1}$ is an automorphism of the Poisson algebra
\[
(\KK[U_+(G_1)], \{.,.\}_{\st}) \otimes  (\KK[U_+(G_2)], \{.,.\}_{\st}) \cong
(\KK[U_+(G_1 \times G_2)], \{.,.\}_{\st}).
\]

(1) In the case when $G_1$ and $G_2$ do not have $SL_2$ factors, \thref{5c} follows at once by applying
\thref{5} to the automorphism $\phi$.

(2) Now consider the case of arbitrary split, connected, simply connected, semisimple algebraic groups $G_1$ and $G_2$. Let
\[
G_i \cong \ol{G}_i \times SL_2^{\times m_i},
\]
where $\ol{G}_i$ are algebraic groups of the same type without $SL_2$ factors and $m_i \in \Nset$. We have the isomorphisms
\[
(\KK[U_+(G_i)], \{.,.\}_{\st}) \cong (\KK[U_+(\ol{G}_i)], \{.,.\}_{\st}) \times (\KK[U_+(SL_2^{m_i})], \{.,.\}_{\st}).
\]

Denote by $\De_+^{i, sl}$ and $\De_+^{i,nsl}$ the subsets of positive simple roots of $G_i$ that 
come from the factors $SL_2^{m_i}$ and $\ol{G}_i$, respectively.
Consider the principal gradings of the Poisson algebras $(\KK[U_+(G_i)], \{.,.\}_{\st})$. 
For the isomorphism in \eqref{eta} and $k \in \Zset_{\geq -1}$, define the linear maps 
\[
\eta_k \colon \KK[U_+(G_1)]_1 \to \KK[U_+(G_2)]_{1+k} \quad
\mbox{so that} \quad \eta|_{\KK[U_+(G_1)]_1} = \sum_{k \geq -1} \eta_k.
\]

The definition of the Poisson structures $\pi_{\st}$ implies that $\KK[U_+(SL_2^{m_i})] \subset \ZZ(\KK[U_+(G_i)])$, where 
as before $\ZZ(.)$ stands for the center of a Poisson algebra. Thus,
\[
Z(\KK[U_+(G_i)]) \cong \ZZ(\KK[U_+(\ol{G}_i)]) \otimes \KK[U_+(SL_2^{m_i})].
\]
For degree reasons, it follows from \prref{c-Poisson-cl} that 
\begin{equation}
\label{1-inc}
\eta_0( x_{\al}) \in \KK[U_+(SL_2^{m_2})]_1, \quad \forall \al \in \De_+^{1, sl}.
\end{equation}
Consider the $\Nset$-gradings of the algebras $(\KK[U_+(G_i)], \{.,.\}_{\st})$ given by 
\[
\deg x_\al =
\begin{cases}
0,  &\mbox{for} \; \; \al \in \De_+^{i, sl}, 
\\
1, &\mbox{for} \; \; \al \in \De_+^{i,nsl}.
\end{cases}
\]
With respect to these gradings $\KK[U_+(G_i)]_0 \subset \ZZ(\KK[U_+(G_i)])$ and by an analogous argument to the one in the proof of \eqref{phi-x}, 
we obtain that
\begin{equation}
\label{2-inc}
\eta_{-1}(x_\al)=0, \; \;  
\eta_0( x_{\al}) \in \KK[U_+(\ol{G}_2)]_1, \quad \forall \al \in \De_+^{1, nsl}.
\end{equation}
It follows from \eqref{1-inc}--\eqref{2-inc} and \prref{gen} that the restrictions of $\eta_0$ to $\KK[U_+(\ol{G}_1)]_1$ and $\KK[U_+(SL_2^{m_1})]_1$
uniquely extend to graded Poisson algebra isomorphisms 
\begin{align*}
&(\KK[U_+(\ol{G}_1)], \{.,.\}_{\st}) \stackrel{\cong}{\lra} (\KK[U_+(\ol{G}_2)], \{.,.\}_{\st}) \quad 
\mbox{and} \\
&(\KK[U_+(SL_2^{m_1})], \{.,.\}_{\st}) \stackrel{\cong}{\lra} (\KK[U_+(SL_2^{m_2})], \{.,.\}_{\st}).
\end{align*}
The second isomorphism implies that $m_1 =m_2$. The first isomorphism and part (1) imply that $\ol{G}_1 \cong \ol{G}_2$. 
Thus, $G_1 \cong G_2$, which completes the proof of the theorem.
\qed
\bre{isom} In a similar fashion one can apply \thref{4} to the isomorphism problem for the class $\CC$ of $\Nset$-graded connected 
cluster algebras of geometric type over a field $\KK$ of characteristic 0, equipped with Gekhtman--Shapiro--Vainshtein Poisson structures.
This class of Poisson algebras is closed under tensor product with respect to the canonical cluster structure on the tensor product 
of cluster algebras. If $\eta \colon (\AA_1, \{.,.\}_1) \stackrel{\cong}{\lra} (\AA_2, \{.,.\}_2)$ is an isomorphism between two Poisson algebras in the class $\CC$, 
then $\phi := \eta \otimes \eta^{-1}$ is an automorphism of their tensor product $(\AA_1 \otimes \AA_2, \{.,.\})$ where $\{.,.\}$ denotes the tensor product
extension of $\{.,.\}_1$ and $\{.,.\}_2$. One can then study the possible forms of $\phi$ by passing to the unipotent automorphisms of 
$(\AA_1 \otimes \AA_2, \{.,.\})$ and applying \thref{4}. 
\ere

\end{document}